\documentclass[10pt]{amsart}

\usepackage{amsmath, amssymb, amsfonts, amsbsy, amsthm, latexsym, graphicx}

\theoremstyle{definition}
\newtheorem{theorem}{Theorem} [section]
\newtheorem{corollary}[theorem]{Corollary}
\newtheorem{lemma}[theorem]{Lemma}
\newtheorem{proposition}[theorem]{Proposition}
\newtheorem{definition}[theorem]{Definition}
\newtheorem{notation}[theorem]{Notation}
\newtheorem{remark}[theorem]{Remark}
\newtheorem{example}[theorem]{Example}
\numberwithin{equation}{section}

\newcommand{\EQ}{\; = \;}
\newcommand{\GE}{\; \ge \;}
\newcommand{\GT}{\; > \;}
\newcommand{\LE}{\; \le \;}
\newcommand{\LT}{\; < \;}
\newcommand{\SUBSET}{\; \subset \;}
\newcommand{\plus}{\; + \;}
\newcommand{\minus}{\; - \;}
\newcommand{\Bc}{{\mathcal{B}}}
\newcommand{\BD}{D_B}
\newcommand{\bigabs}[1]{{\bigl|#1\bigr|}}
\newcommand{\Bigabs}[1]{{\Bigl|#1\Bigr|}}
\newcommand{\biggabs}[1]{{\biggl|#1\biggr|}}
\newcommand{\biggbracket}[1]{\biggl[#1\biggr]}
\newcommand{\C}{\mathbf{C}}
\newcommand{\clspan}{{\overline{\mbox{\rm span}}}}
\newcommand{\comp}{{\mathrm{C}}}
\newcommand{\Dim}{{\mathrm{dim}}}
\newcommand{\dist}{{\mathrm{dist}}}
\newcommand{\Ec}{{\mathcal{E}}}
\newcommand{\tEc}{{\tilde{\mathcal{E}}}}
\newcommand{\te}{{\tilde{e}}}
\newcommand{\eps}{\varepsilon}
\newcommand{\Fc}{{\mathcal{F}}}
\newcommand{\tFc}{{\tilde{\mathcal{F}}}}
\newcommand{\tf}{{\tilde{f}}}
\newcommand{\Gc}{{\mathcal{G}}}
\newcommand{\tg}{{\tilde{g}}}
\newcommand{\ip}[2]{\langle#1,#2\rangle}
\newcommand{\bigip}[2]{\bigl\langle #1, \, #2 \bigr\rangle}
\newcommand{\biggip}[2]{\biggl\langle #1, \, #2 \biggr\rangle}
\newcommand{\Int}{{\mathrm{Int}}}
\newcommand{\cM}{{\mathcal{M}}}
\newcommand{\one}{\mathbf{1}}
\newcommand{\N}{\mathbf{N}}
\newcommand{\norm}[1]{\|#1\|}
\newcommand{\bignorm}[1]{\bigl\|#1\bigr\|}
\newcommand{\Bignorm}[1]{\Bigl\|#1\Bigr\|}
\newcommand{\biggnorm}[1]{\biggl\|#1\biggr\|}
\newcommand{\Pb}{\mathbf{P}}
\newcommand{\bigparen}[1]{\bigl(#1\bigr)}
\newcommand{\Bigparen}[1]{\Bigl(#1\Bigr)}
\newcommand{\biggparen}[1]{\biggl(#1\biggr)}
\newcommand{\plim}{\operatornamewithlimits{\mbox{$p$}-\mathrm{lim}}}

\newcommand{\Qc}{{\mathcal{Q}}}
\newcommand{\R}{\mathbf{R}}
\newcommand{\rank}{{\mathrm{rank}}}
\newcommand{\set}[1]{\{#1\}}
\newcommand{\bigset}[1]{\bigl\{#1\bigr\}}
\newcommand{\Bigset}[1]{\Bigl\{#1\Bigr\}}
\newcommand{\Span}{\mathrm{span}}
\newcommand{\spectrum}{\mathrm{Sp}}

\newcommand{\trace}{{\mathrm{trace}}}
\newcommand{\Z}{\mathbf{Z}}

\begin{document}

\title[Density, Overcompleteness, and Localization, I]
{Density, Overcompleteness, and \\ Localization of Frames. \\ I. Theory}

\author[R.~Balan, P.~G.~Casazza, C.~Heil, and Z.~Landau]
{Radu~Balan, Peter~G.~Casazza, Christopher~Heil, and Zeph~Landau}

\address{\textrm{(R.~Balan)}
Siemens Corporate Research, 
755 College Road East, 
Princeton, NJ 08540}
\email{radu.balan@siemens.com}

\address{\textrm{(P.~G.~Casazza)}
Department of Mathematics,
University of Missouri,
Columbia, MO 65211}
\email{pete@math.missouri.edu}

\address{\textrm{(C.~Heil)}
School of Mathematics,
Georgia Institute of Technology,
Atlanta, GA 30332}
\email{heil@math.gatech.edu}

\address{\textrm{(Z.~Landau)}
Department of Mathematics R8133,
The City College of New York,
Convent Ave at 138th Street,
New York, NY 10031}
\email{landau@sci.ccny.cuny.edu}

\date{March 30, 2005}

\keywords{
Density, excess, frames, Gabor systems, modulation spaces, overcompleteness,
Riesz bases, wavelets, Weyl--Heisenberg systems.
}

\subjclass[2000]{Primary 42C15; Secondary 46C99}

\thanks{
The second author was partially supported by NSF Grants
DMS-0102686 and DMS-0405376.
The third author was partially supported by NSF Grant DMS-0139261.
Some of the results of this paper were previously announced, without proofs,
in the research announcement \cite{BCHL05b}.}

\begin{abstract}
This work presents a quantitative framework for describing the
overcompleteness of a large class of frames.
It introduces notions of localization and approximation between two frames
$\mathcal{F} = \{f_i\}_{i \in I}$ and
$\mathcal{E} = \{e_j\}_{j \in G}$ ($G$ a discrete abelian group),
relating the decay of the expansion of the elements of~$\mathcal{F}$
in terms of the elements of $\mathcal{E}$ via a map $a \colon I \to G$. 
A fundamental set of equalities are shown between three seemingly unrelated
quantities: the relative measure of $\mathcal{F}$,
the relative measure of $\mathcal{E}$ ---
both of which are determined by certain averages of inner products
of frame elements with their corresponding dual frame elements ---
and the density of the set $a(I)$ in $G$.
Fundamental new results are obtained on the excess and overcompleteness
of frames, on the relationship between frame bounds and density,
and on the structure of the dual frame of a localized frame.
In a subsequent paper, these results are applied to the case of
Gabor frames, producing an array of new results as well as clarifying
the meaning of existing results.

The notion of localization and related approximation properties
introduced in this paper are a spectrum of ideas that quantify the
degree to which elements of one frame can be approximated by elements of
another frame.
A comprehensive examination of the interrelations among these localization
and approximation concepts is presented.
\end{abstract}

\copyrightinfo{}{}

\maketitle

\section{Introduction}

The fundamental structural feature of frames that are not Riesz bases is
the overcompleteness of its elements.
To date, even partial understanding of this overcompleteness has been
restricted to limited examples, such as finite-dimensional frames, frames
of windowed exponentials, or frames of time-frequency shifts (Gabor systems).
The ideas and results presented here provide a quantitative framework for
describing the overcompleteness of a large class of frames.
The consequences of these ideas are:
(a)~an array of fundamental new results for frames that hold in a
general setting,
(b)~significant new results for the case of Gabor frames,
as well as a new framing of existing results that clarifies their meaning,
and (c)~the presentation of a novel and fruitful point of view for
future research.

Due to the length of this work, it is natural to present it in two parts.
The first part, containing the theoretical and structural results that have
driven the research, forms this paper.
The second part, containing the applications to Gabor frames,
will appear in the paper \cite{BCHL05a} (hereafter referred to as ``Part~II'').

At the core of our main results is Theorem~\ref{densityredundancy}.
The precise statement of the theorem requires some detailed notation, but
the essence of the result can be summarized as follows.
We begin with two frames $\Fc = \set{f_i}_{i \in I}$ and
$\Ec = \set{e_j}_{j \in G}$, where $G$ is a discrete abelian group,
and introduce a notion of the localization of~$\Fc$ with respect to $\Ec$.
The idea of localization is that it describes the decay of the coefficients
of the expansion of elements of $\Fc$ in terms of the elements of $\Ec$.
To make this notion of decay meaningful, a map $a$ from the index set $I$ into
the index set $G$ is introduced.
With this setup, Theorem~\ref{densityredundancy} establishes a remarkable
equality relating three seemingly unrelated quantities: certain averages
of $\ip{f_i}{\tf_i}$ and $\ip{e_j}{\te_j}$
of frame elements with corresponding canonical dual frame elements,
which we refer to as \emph{relative measures},
and the density of the set $a(I)$ in $G$.
This equality between density and relative measure is striking since
the relative measure is a function of the frame elements,
while the density is solely determined by the index set~$I$ and
the mapping $a \colon I \to G$.

The impact of Theorem~\ref{densityredundancy} comes in several forms.
First, the result itself is new, and its consequences along with related
ideas discussed in more detail below represent a significant increase in
the understanding of the structure of abstract frames.
Second, the application of Theorem~\ref{densityredundancy} and our other new
theorems to the case of Gabor frames yields new results,
which will be presented in Part~II.
These recover as corollaries the existing density results known to hold
for Gabor frames, but in doing so, shows them in a new light,
as the consequence of more general considerations rather than of a
particular rigid structure of the frames themselves.
The notions of localization, approximation, and measure are
interesting and useful new ideas which we feel will have impact beyond
the results presented in this paper.
In particular, it will be interesting to see to what degree wavelet frames
fit into this framework, especially given recent results on density theorems
for affine frames \cite{HK03}, \cite{SZ02}.

In addition to the fundamental equalities relating density and measures
discussed above, we obtain a set of additional significant results,
as follows.

First, we provide a comprehensive theory of \emph{localization} of frames.
Localization is not a single concept, but a suite of related ideas.
We introduce a collection of definitions and describe the implications among
these various definitions.
We also introduce a set of \emph{approximation properties} for frames,
and analyze the interrelations between these properties and the localization
properties.

Second, we explore the implications of the connection between density and
overcompleteness.
We show that in any overcomplete frame which possesses sufficient localization,
the overcompleteness must have a certain degree of uniformity.
Specifically, we construct an infinite subset of the frame with positive
density which can be removed yet still leave a frame.
We obtain relations among the frame bounds, density of the index set $I$,
and norms of the frame elements, and prove in particular that if $\Fc$ is
a tight localized frame whose elements all have the same norm then
the index set $I$ must have uniform density.

Third, we explore the structure of the dual frame, showing that if a frame
is sufficiently localized then its dual frame is also.
We also prove that any sufficiently localized frame can be written as a finite
union of Riesz sequences.
This shows that the Feichtinger conjecture
(which has recently been shown to be equivalent to the
famous \emph{Kadison--Singer conjecture} \cite{CT05})
is true for the case of localized frames.

In Part~II we apply our results to derive new implications for the case of
Gabor frames and more general systems of Gabor molecules, whose elements
are not not required to be simple time-frequency shifts of each other,
but instead need only share a common envelope of concentration about points
in the time-frequency plane.
These include strong results on the the structure of the dual frame
of an irregular Gabor frame, about which essentially nothing has
previously been known beyond the fact that it consists of a set of
$L^2$ functions.
We prove that if an irregular Gabor frame is generated by a function $g$
which is sufficiently concentrated in the time-frequency plane
(specifically, $g$ lies in the modulation space $M^1$),
then the elements of the dual frame also lie in $M^1$.
We further prove that the dual frame forms a set of Gabor molecules,
and thus, while it need not form a Gabor frame, the elements do share
a common envelope of concentration in the time-frequency plane.
Moreover, this same result applies if the original frame was only
itself a frame of Gabor molecules.

Our paper is organized as follows.
The next subsection will give a more detailed and precise summary and outline
of our results.
Section~\ref{section2} introduces the concepts of localization and
approximation properties and presents the interrelations among them.
We also define density and relative measure precisely in
that section.
The main results of this paper for abstract frames are presented in
Section~\ref{mainsec}.

\subsection{Outline} \label{outline}

\subsubsection{Density, Localization, HAP, and Relative Measure}

The main body of our paper begins in Section~\ref{section2}, where,
following the definition of density in Section~\ref{densitysec}, we define
several types of localization and approximation properties for abstract frames
in Sections~\ref{localizationsec} and~\ref{approxsec}.

Localization is determined both by the frame $\Fc = \set{f_i}_{i \in I}$
and by a reference system $\Ec = \set{e_j}_{j \in G}$.
We assume the reference system is indexed by a group of the form
\begin{equation} \label{group}
G \EQ \prod_{i=1}^d a_i \Z \, \times \, \prod_{j=1}^e \Z_{n_j},
\end{equation}
with a metric on $G$ defined as follows.
If $m_j \in \Z_{n_j}$, set $\delta(m_j) = 0$ if $m_j =0$, otherwise
$\delta(m_j) = 1$.
Then given $g = (a_1 n_1, \dots, a_d n_d, m_1, \dots, m_e) \in G$, set
\begin{equation} \label{metric}
|g|
\EQ \sup\bigset{|a_1n_1|, \, \dots,\, |a_d n_d|, \,
                \delta(m_1), \, \dots, \, \delta(m_e)}.
\end{equation}
The metric is then $d(g,h) = |g-h|$ for $g$, $h \in G$.
Our results can be generalized to other groups;
the main properties of the group defined by \eqref{group} that are used are
that $G$ is a countably infinite abelian group which has a shift-invariant
metric with respect to which it is locally finite.
The reader can simply take $G = \Z^d$ without much loss of insight
on a first reading.

The additive structure of the index set $G$ of the reference system does play
a role in certain of our results.
However, the index set $I$ of the frame $\Fc$ need not be structured.
For example, in our applications in Part~II we will have an irregular
Gabor system
$\Fc = \Gc(g,\Lambda)
= \set{e^{2\pi i \eta x} g(x-u)}_{(u,\eta) \in \Lambda}$,
which has as its index set an arbitrary countable subset
$\Lambda \subset \R^{2d}$,
while our reference system will be a lattice Gabor system
$\Ec = \Gc(\phi,\alpha\Z^d \times \beta\Z^d)
= \set{e^{2\pi i \eta x} \phi(x-u)}
  _{(u,\eta) \in \alpha\Z^d \times \beta\Z^d}$,
indexed by $G = \alpha\Z^d \times \beta\Z^d$.

A set of approximation properties for abstract frames is introduced
in Definition~\ref{approxdef}.
These are defined in terms of how well the elements
of the reference system are approximated by finite linear combinations
of frame elements, or vice versa, and provide an abstraction for
general frames of the essential features of the Homogeneous Approximation
Property (HAP) that is known to hold for Gabor frames or windowed exponentials
(see \cite{RS95}, \cite{GR96}, \cite{CDH99}).

We list in Theorem~\ref{relations} the implications that hold among the
localization and approximation properties.
In particular, there is an equivalence between
$\ell^2$-column decay and the HAP, and between
$\ell^2$-row decay and a dual HAP.

In Section~\ref{selflocsec} we introduce another type of localization.
Instead of considering localization with respect to a fixed reference sequence,
we consider localizations in which the reference is the frame itself
(``self-localization'') or its own canonical dual frame.
Theorem~\ref{selflocthm} states that every $\ell^1$-self-localized frame
is $\ell^1$-localized with respect to its canonical dual frame.
The proof of this result is an application of a type of noncommutative
Wiener's Lemma, and is given in Appendix~\ref{selflocappend}.

We define the density of an abstract frame $\Fc = \set{f_i}_{i \in I}$
in Section~\ref{densitysec}.
We assume there is some associated mapping $a \colon I \to G$.
For example, in the Gabor case, $I = \Lambda$ is an arbitrary countable
sequence in $\R^{2d}$ while $G = \alpha\Z^d \times \beta\Z^d$,
and $a$ maps elements of $I$ to elements of~$G$ by rounding off to a
near element of $G$ (note that $a$ will often not be injective).
Then density is defined by considering the average number of points in $a(I)$
inside boxes of larger and larger radius.
By taking the infimum or supremum over all boxes of a given radius and then
letting the radius increase, we obtain lower and upper densities
$D^\pm(I,a)$.
By using limits with respect to an ultrafilter $p$ 
and a particular choice of centers $c = \set{c_N}_{N \in \N}$ for the boxes,
we obtain an entire collection of densities
$D(p,c)$ intermediate between the upper and lower densities
(for background on ultrafilters, we refer to \cite[Chap.~3]{HS98}
or \cite[App.~A]{BCHL05a}).

The relative measure of an abstract frame sequence $\Fc = \set{f_i}_{i \in I}$
with respect to a reference frame sequence $\Ec = \set{e_j}_{j \in G}$
is introduced in Section~\ref{relativemeasure}.
For simplicity, in this introduction we discuss only the case where
both are frames for the entire space; in this case we speak of the
\emph{measures} of $\Fc$ and $\Ec$ instead of the \emph{relative measures}.
Furthermore we will discuss here only the case where $\Ec$ is a Riesz basis,
so that its measure is $1$.
Let $S_N(j)$ denote the discrete ``box'' in $G$ centered at $j \in G$
and with ``side lengths'' $N$ (see equation \eqref{boxdef} for the
precise definition).
Let $I_N(j) = a^{-1}(S_N(j))$ denote the preimage in $I$ of $S_N(j)$
under the map~$a \colon I \to G$.
We declare the lower measure of the frame $\Fc$ to be
$$\cM^-(\Fc)
\EQ \liminf_{N \to \infty} \, \inf_{j \in G} \,
    \frac1{|I_N(j)|} \sum_{i \in I_N(j)} \ip{f_i}{\tf_i},$$
and make a similar definition for the upper measure $\cM^+(\Fc)$
(note that $0 \le \ip{f_i}{\tf_i} \le 1$ for all~$i$).
We also define the measure $\cM(\Fc;p,c)$ with respect to an ultrafilter $p$
and a particular choice of box centers $c = (c_N)_{N \in \N}$.
Thus, as was the case with the densities, we actually have a suite of
definitions, a range of measures that are intermediate between
the lower and upper measures.
Note that if $\Fc$ is a Riesz basis, then $\ip{f_i}{\tf_i} = 1$ for
every~$i$, so a Riesz basis has upper and lower measure~$1$.
The definition of relative measure becomes more involved when the systems
are only frame sequences, i.e., frames for their closed linear spans.
In this case, the relative measures are determined by averages of
$\ip{P_\Ec f_i}{\tf_i}$ or $\ip{P_\Fc \te_j}{e_j}$, respectively,
where $P_\Ec$ and $P_\Fc$ are the orthogonal projections onto the
closed spans of $\Ec$ and $\Fc$.
The precise definition is given in Definition~\ref{redundancydef}.

\subsubsection{Density and Overcompleteness for Localized Frames}
Section~\ref{necessarysec} presents two necessary conditions on
the density of a frame.
In Theorem~\ref{necessary}, we show that a frame which satisfies only a weak 
version of the HAP will satisfy a Nyquist-type condition, specifically,
it must have a lower density which satisfies $D^-(I,a) \ge 1$.
In Theorem~\ref{finitedensity}, we show that under a stronger localization
assumption, the upper density must be finite.

The connection between density and overcompleteness, which is among the most
fundamental of our main results, is presented in Section~\ref{connectsec}.
We establish a set of equalities between the relative measures
and the reciprocals of the density.
Specifically, we prove in Theorem~\ref{densityredundancy} that for frame
$\Fc$ that is appropriately localized with respect to a Riesz basis~$\Ec$,
we have the following equalities for the lower and upper measures and
for every measure defined with respect to an ultrafilter $p$
and sequence of centers $c = (c_N)_{N \in \N}$ in $G$:
\begin{equation} \label{equalities}
\cM^-(\Fc) \EQ \frac1{D^+(I,a)}, \quad
\cM(\Fc;p,c) \EQ \frac1{D(p,c)}, \quad
\cM^+(\Fc) \EQ \frac1{D^-(I,a)}.
\end{equation}
Moreover, we actually obtain much finer versions of the equalities above which
hold for the case of a frame sequence compared to a reference system that is
also a frame sequence.
The left-hand side of each equality is a function of the frame elements,
while the right-hand side is determined by the index set alone.
As immediate consequences of these equalities we obtain inequalities relating
density, frame bounds, and norms of the frame elements.
In particular, we show that if $\Fc$ and $\Ec$ are both localized
tight frames whose frame elements all have identical norms, then the
index set $I$ must have uniform density, i.e., the upper and lower
densities of $I$ must be equal.
Thus tightness necessarily requires a certain uniformity of the index set.

The equalities in \eqref{equalities} suggest that relative measure is a
quantification of overcompleteness for localized frames.
To illustrate this connection, let us recall the definition of the
\emph{excess} of a frame, which is a crude measure of overcompleteness.
The excess of a frame $\set{f_i}_{i \in I}$ is the cardinality of the
largest set $J$ such that $\set{f_i}_{i \in I \setminus J}$ is complete
(but not necessarily still a frame).
An earlier paper \cite{BCHL03} showed that there is an infinite
$J \subset I$ such that $\set{f_i}_{i \in I \setminus J}$ is still a frame
if and only if there exists an infinite set $J_0 \subset I$ such that
\begin{equation} \label{excesscondition}
\sup_{i \in J_0} \, \ip{f_i}{\tf_i} \LT 1.
\end{equation}
The set $J$ to be removed will be a subset of $J_0$, but, in general,
the technique of \cite{BCHL03} will construct only an extremely sparse
set $J$ (typically zero density in the terminology of this paper).
If $\cM^-(\Fc) < 1$, then \eqref{excesscondition} will be satisfied for
some $J_0$ (see Proposition~\ref{infiniteexcess}), and so some infinite
set can be removed from the frame.
We prove in Section~\ref{positivesec} that if a frame is localized and
$\cM^+(\Fc) < 1$, then not merely can some infinite set be removed,
but this set can be chosen to have positive density.
We believe, although we cannot yet prove, that the reciprocal of the
relative measure is in fact quantifying the redundancy of an abstract frame,
in the sense that it should be the case that if $\Fc$ is appropriately
localized and $\cM^+(\Fc) < 1$, then there should be a subset of $\Fc$
with density $\frac1{\cM^+(\Fc)} - 1 - \eps$ which can be removed leaving
a subset of $\Fc$ with density $1+\eps$ which is still a frame for~$H$.

The last of our results deals with the conjecture of Feichtinger
that every frame that is norm-bounded below
($\inf_i \norm{f_i} > 0$) can be written as a union of a finite number of
Riesz sequences (systems that are Riesz bases for their closed linear spans).
It is shown in \cite{CCLV03}, \cite{CV03}, \cite{CT05}
that Feichtinger's conjecture equivalent to the
celebrated Kadison--Singer (paving) conjecture.
In Section~\ref{rieszsec}, we prove that this conjecture is true for
the case of $\ell^1$-self-localized frames which are norm-bounded below.
This result is inspired by a similar result of Gr\"ochenig's from
\cite{Gro03} for frames which are sufficiently localized in his sense,
although our result is distinct.
Another related recent result appears in \cite{BS04}.

We believe that localization is a powerful and useful new concept.
As evidence of this fact, we note that Gr\"ochenig has independently
introduced a concept of localized frames, for a completely different
purpose \cite{Gro04}.
We learned of Gr\"ochenig's results shortly after completion of our own
major results.
The definitions of localizations presented here and in \cite{Gro04}
differ, but the fact that this single concept has independently arisen
for two very distinct applications shows its utility.
In his elegant paper, Gr\"ochenig has shown that frames which are sufficiently
localized in his sense provide frame expansions not only for the Hilbert
space~$H$ but for an entire family of associated Banach function spaces.
Gr\"ochenig further showed that if a frame is sufficiently localized
in his sense (a~polynomial or exponential localization) then the dual frame
is similarly localized.

\subsection{General Notation} \label{prelimsec}

The following notation will be employed throughout this paper.
$H$ will refer to a separable Hilbert space,
$I$ will be a countable index set, and~$G$ will be the group given
by \eqref{group} with the metric defined in \eqref{metric}.
We implicitly assume that there exists a map $a \colon I \to G$
associated with $I$ and~$G$.
The map $a$ induces a semi-metric
$d(i,j) = |a(i) - a(j)|$ on $I$.
This is only a semi-metric since $d(i,j) = 0$ need not imply $i=j$.

The finite linear span of a subset $S \subset H$ is denoted $\Span(S)$,
and the closure of this set is $\clspan(S)$.
The cardinality of a finite set $E$ is denoted by $|E|$.

For each integer $N > 0$ we let
\begin{equation} \label{boxdef}
S_N(j)
\EQ \Bigset{k \in G : |k-j| \le \frac{N}2}
\end{equation}
denote a discrete ``cube'' or ``box'' in $G$ centered at $j \in G$.
The cardinality of $S_N(j)$ is independent of $j$.
For example, if $G = \Z^d$ then
$|S_{2N}(j)| = |S_{2N+1}(j)| = (2N+1)^d$.
In general, there will exist a constant $C$ and integer $d>0$ such that
\begin{equation} \label{asymptotics}
\lim_{N \to \infty} \frac{|S_N(j)|}{N^d} \EQ C.
\end{equation}
We let $I_N(j)$ denote the inverse image of $S_N(j)$ under $a$, i.e.,
$$I_N(j) \EQ a^{-1}(S_N(j))
\EQ \set{i \in I : a(i) \in S_N(j)}.$$

\subsection{Notation for Frames and Riesz Bases}

We use standard notations for frames and Riesz bases as found in the texts
\cite{Chr03}, \cite{Dau92}, \cite{Gro01}, \cite{You01}
or the research-tutorials \cite{Cas00}, \cite{HW89}.
Some particular notation and results that we will need are as follows.

A sequence $\Fc = \set{f_i}_{i \in I}$ is a \emph{frame} for $H$ if there
exist constants $A$, $B > 0$, called \emph{frame bounds}, such that
\begin{equation} \label{framedef}
\forall\, f \in H, \quad
A \, \norm{f}^2 \LE \sum_{i \in I} |\ip{f}{f_i}|^2 \LE B \, \norm{f}^2.
\end{equation}
The \emph{analysis operator} $T \colon H \to \ell^2(I)$ is
$Tf = \set{\ip{f}{f_i}}_{i \in I}$, and its adjoint
$T^* c = \sum_{i \in I} c_i \, f_i$
is the \emph{synthesis operator}.
The \emph{Gram matrix} is
$T T^* = [\ip{f_i}{f_j}]_{i,j \in I}$.
The \emph{frame operator}
$Sf = T^* T f = \sum_{i \in I} \ip{f}{f_i} \, f_i$
is a bounded, positive, and invertible mapping of $H$ onto itself.
The \emph{canonical dual frame} is
$\tFc = S^{-1}(\Fc) = \set{\tf_i}_{i \in I}$
where $\tf_i = S^{-1} f_i$.
For each $f \in H$ we have the \emph{frame expansions}
$f = \sum_{i \in I} \ip{f}{f_i} \, \tf_i
= \sum_{i \in I} \ip{f}{\tf_i} \, f_i$.
We call $\Fc$ a \emph{tight frame} if we can take $A=B$, and
a \emph{Parseval frame} if we can take $A=B=1$.
If $\Fc$ is any frame, then $S^{-1/2}(\Fc)$ is the
\emph{canonical Parseval frame} associated to $\Fc$.
We call $\Fc$ a \emph{uniform norm frame} if all the frame elements have
identical norms, i.e., if $\norm{f_i} = const.$ for all $i \in I$.

A sequence which satisfies the upper frame bound estimate
in \eqref{framedef}, but not necessarily the lower estimate, is called
a \emph{Bessel sequence} and $B$ is a \emph{Bessel bound}.
In this case, $\norm{\sum c_i f_i}^2 \le B \, \sum |c_i|^2$ for any
$(c_i)_{i \in I} \in \ell^2(I)$.
In particular, $\norm{f_i}^2 \le B$ for every $i \in I$, i.e., all
Bessel sequences are norm-bounded above.
If we also have $\inf_i \norm{f_i} > 0$, then we say the sequence is
norm-bounded below.

We will also consider sequences that are frames for their closed linear spans
instead of for all of $H$.
We refer to such a sequence as a \emph{frame sequence}.
If $\Fc = \set{f_i}_{i \in I}$ is a frame sequence,
then $\tFc = \set{\tf_i}_{i \in I}$ will denote its canonical dual frame
within $\clspan(F)$.
The orthogonal projection $P_\Fc$ of $H$ onto $\clspan(\Fc)$ is given~by
\begin{equation} \label{orthogproj}
P_\Fc f \EQ \sum_{i \in I} \ip{f}{f_i} \, \tf_i,
\qquad f \in H.
\end{equation}

A frame is a basis if and only if it is a Riesz basis, i.e., the image
of an orthonormal basis for $H$ under a continuous, invertible linear mapping.
We say $\Fc = \set{f_i}_{i \in I}$ is a \emph{Riesz sequence} if it is a
Riesz basis for its closed linear span in~$H$.
In this case the canonical dual frame $\tFc = \set{\tf_i}_{i \in I}$
is the unique sequence in $\clspan(\Fc)$ that is biorthogonal to~$\tFc$, i.e.,
$\ip{f_i}{\tf_j} = \delta_{ij}$.

\section{Density, Localization, HAP, and Relative Measure} \label{section2}

\subsection{Density} \label{densitysec}

Given an index set $I$ and a map $a \colon I \to G$, we define
the density of~$I$ by computing the analogue of Beurling density of its
image $a(I)$ as a subset of~$G$.
Note that we regard $I$ as a sequence, and hence repetitions of images
count in determining the density.
The precise definition is as follows.

\begin{definition}[Density] \label{densitydef}
The \emph{lower and upper densities of $I$ with respect to $a$} are
\begin{equation} \label{lowerdensity}
D^-(I,a)
\EQ \liminf_{N \to \infty} \inf_{j \in G} \frac{|I_N(j)|}{|S_N(j)|},
\qquad
D^+(I,a)
\EQ \limsup_{N \to \infty} \sup_{j \in G} \frac{|I_N(j)|}{|S_N(j)|},
\end{equation}
respectively.
Note that these quantities could be zero or infinite, i.e., we have
$0 \le D^-(I,a) \le D^+(I,a) \le \infty$.
When $D^-(I,a) = D^+(I,a) = D$ we say $I$ has \emph{uniform density}~$D$.
~\qed
\end{definition}

These lower and upper densities are only the extremes of the
possible densities that we could naturally assign to $I$ with respect to $a$.
In particular, instead of taking the infimum or supremum over all possible
centers in \eqref{lowerdensity} we could choose one
specific sequence of centers, and instead of computing the liminf or limsup
we could consider the limit with respect to some ultrafilter.
The different possible choices of ultrafilters and sequences of centers
gives us a natural collection of definitions of density, made precise in
the following definition.

\begin{definition}
Let $p$ be a free ultrafilter, and let $c = (c_N)_{N \in \N}$ be any
sequence of centers $c_N \in G$.
Then the \emph{density of $I$ with respect to $a$, $p$, and $c$} is
$$D(p,c)
\EQ D(p,c;I,a)
\EQ \plim_{N \in \N} \frac{|I_N(c_N)|}{|S_N(c_N)|}.
\quad\qed$$
\end{definition}

\begin{example}
If $I=G$ and $a$ is the identity map, then $I_N(j) = S_N(j)$ for every $N$
and~$j$, and hence $D(p,c) = D^-(I,a) = D^+(I,a) = 1$
for every choice of free ultrafilter $p$ and sequence of centers~$c$.
\qed
\end{example}

The following example shows how the density we have defined relates to the
standard Beurling density of the index set of a Gabor system.

\begin{example}[Gabor Systems] \label{gabordensityrel}
Consider a Gabor system $\Fc = \Gc(g,\Lambda)$ and a reference Gabor system
$\Ec = \Gc(\phi,\alpha\Z^d \times \beta\Z^d)$.
The index set $I = \Lambda$ is a countable sequence of points in $\R^{2d}$,
and the reference group is $G = \alpha\Z^d \times \beta\Z^d$.
A natural map $a \colon \Lambda \to G$ is a
simple roundoff to a near element of $G$, i.e.,
$$a(x,\omega)
\EQ \bigparen{\alpha \, \Int\bigparen{\tfrac{x}{\alpha}},
              \beta \, \Int\bigparen{\tfrac{\omega}{\beta}}},
\qquad (x,\omega) \in \Lambda,$$
where $\Int(x) = (\lfloor x_1 \rfloor, \dots, \lfloor x_d \rfloor)$.
With this setup, $S_N(j)$ is the intersection of
$\alpha\Z^d \times \beta\Z^d$ with the cube $Q_N(j)$ in $\R^{2d}$
centered at $j$ with side lengths $N$.
Such a cube contains approximately $N^{2d}/(\alpha\beta)^d$ points of
$\alpha\Z^d \times \beta\Z^d$; precisely,
$$\lim_{N \to \infty} \frac{|S_N(j)|}{N^{2d}} \EQ \frac1{(\alpha\beta)^d}.$$
Also, because $a$ is a bounded perturbation of the identity map,
the number of points in $I_N(j)$ is asymptotically the cardinality of
$\Lambda \cap Q_N(j)$.
Consequently, the standard definition of the upper Beurling density
$\BD^+(\Lambda)$ of $\Lambda$ is related to our definition of the
upper density of $\Lambda$ with respect to $a$ as follows:
\begin{align*}
\BD^+(\Lambda)
& \EQ \limsup_{N \to \infty} \sup_{j \in \R^{2d}}
      \frac{|\Lambda \cap Q_N(j)|}{N^{2d}} \\
& \EQ \frac1{(\alpha\beta)^d} \, \limsup_{N \to \infty}
      \sup_{j \in \alpha\Z^d \times \beta\Z^d}
      \frac{|I_N(j)|}{|S_N(j)|}
\EQ \frac1{(\alpha\beta)^d} \, D^+(\Lambda,a). \notag
\end{align*}
Similarly the lower Beurling density of $\Lambda$ is
$\BD^-(\Lambda) = (\alpha\beta)^{-d} \, D^-(\Lambda,a)$.
In particular, when $\alpha\beta = 1$ (the ``critical density'' case),
our definition coincides with Beurling density, but in general the extra
factor of $(\alpha\beta)^d$ must be taken into account.
\qed
\end{example}

The following two lemmas will be useful later for our density calculations.
The first lemma is similar to \cite[Lem.~20.11]{HS98}.

\begin{lemma} \label{liminfexist}
Let $a \colon I \to G$ be given.

\smallskip
\begin{enumerate}
\item[(a)]
For every free ultrafilter $p$ and sequence of centers
$c = (c_N)_{N \in \N}$ in $G$, we have
$D^-(I,a) \le D(p,c) \le D^+(I,a)$.

\smallskip
\item[(b)]
There exist free ultrafilters $p^-$, $p^+$ and sequence of centers
$c^- = (c_N^-)_{N \in \N}$, $c^+ = (c_N^+)_{N \in \N}$ in $G$ such that
$D^-(I,a) = D(p^-,c^-)$
and
$D^+(I,a) = D(p^+,c^+)$.
\end{enumerate}
\end{lemma}
\begin{proof}
(a) Follows immediately from the properties of ultrafilters.

\medskip
(b) For each $N > 0$, we can choose a point $c_N$ so that
$$\inf_{j \in G} \frac{|I_N(j)|}{|S_N(j)|}
\LE \frac{|I_N(c_N)|}{|S_N(j)|}
\LE \Bigparen{\inf_{j \in G} \frac{|I_N(j)|}{|S_N(j)|}} \plus \frac1N.$$
Then we can choose a free ultrafilter $p$ such that
$$\plim_{N \in \N} \frac{|I_N(c_N)|}{|S_N(j)|}
\EQ \liminf_{N \to \infty} \frac{|I_N(c_N)|}{|S_N(j)|}.$$
For these choices, we have
\begin{align*}
D^-(I,a) \LE D(p,c)
& \EQ \plim_{N \in \N} \frac{|I_N(c_N)|}{|S_N(j)|} \\[1 \jot]
& \LE \liminf_{N \to \infty} \, \biggbracket{
      \Bigparen{\inf_{j \in G} \frac{|I_N(j)|}{|S_N(j)|}} \plus \frac1N}
      \\[1 \jot]
& \LE \Bigparen{\liminf_{N \to \infty} \inf_{j \in G}
      \frac{|I_N(j)| }{|S_N(j)|}}
      \plus \Bigparen{\limsup_{N \to \infty} \frac1N}
\EQ D^-(I,a).
\end{align*}
Thus we can take $p^- = p$ and $c^- = (c_N)_{N \in \N}$.
The construction of $p^+$ and $c^+$ is similar.
\end{proof}

\begin{lemma} \label{annuluslemma}
Assume $D^+(I,a) < \infty$.
Then $K = \sup_{j \in G} |a^{-1}(j)|$ is finite, and
for any set $E \subset G$ we have
\begin{equation} \label{annulus}
|a^{-1}(E)| \LE K \, |E|.
\end{equation}
\end{lemma}

\subsection{The Localization Properties} \label{localizationsec}

We now introduce a collection of definitions of localization,
given in terms of the decay of the inner products of the elements of one
sequence $\Fc$ with respect to the elements of a reference sequence $\Ec$.
In Section~\ref{approxsec}, we define several approximation properties,
which are determined by how well the elements of one sequence are approximated
by finite linear combinations of the elements of the other sequence.
The relationships among these properties is stated in Theorem~\ref{relations}.

The words ``column'' and ``row'' in the following definition
refer to the $I \times G$ cross-Grammian matrix
$[\ip{f_i}{e_j}]_{i \in I, j \in G}$.
We think of the elements in locations $(i,a(i))$ as corresponding to the
main diagonal of this matrix.

\begin{definition}[Localization] \label{localizationdef}
Let $\Fc = \set{f_i}_{i \in I}$ and
$\Ec = \set{e_j}_{j \in G}$ be sequences in~$H$
and $a \colon I \to G$ an associated map.

\smallskip
\begin{enumerate}
\item[(a)]
We say $\Fc$ is \emph{$\ell^p$-localized} with respect to the
reference sequence $\Ec$ and the map $a$,
or simply that \emph{$(\Fc,a,\Ec)$ is $\ell^p$-localized}, if
$$\sum_{j \in G} \, \sup_{i \in I} \, |\ip{f_i}{e_{j+a(i)}}|^p
\LT \infty.$$
Equivalently, there must exist an $r \in \ell^p(G)$ such that
$$\forall\, i \in I, \quad
\forall\, j \in G, \quad
|\ip{f_i}{e_j}| \LE r_{a(i)-j}.$$

\medskip
\item[(b)]
We say that $(\Fc,a,\Ec)$ has \emph{$\ell^p$-column decay} if for every
$\eps > 0$ there is an integer $N_\eps > 0$ so that
\begin{equation} \label{PropXdef}
\forall\, j \in G, \quad
\sum_{i \in I \setminus I_{N_\eps}(j)} |\ip{f_i}{e_j}|^p \LT \eps.
\end{equation}

\medskip
\item[(c)]
We say $(\Fc,a,\Ec)$ has \emph{$\ell^p$-row decay} if for every
$\eps > 0$ there is an integer $N_\eps > 0$ so that
\begin{equation} \label{PropXstardef}
\forall\, i \in I, \quad
\sum_{j \in G \setminus S_{N_\eps}({a(i)})} |\ip{f_i}{e_j}|^p \LT \eps.
\quad\qed
\end{equation}
\end{enumerate}
\end{definition}

Note that given a sequence $\Fc$, the definition of localization is dependent 
upon both the choice of reference sequence $\Ec$ and the map $a$.

\begin{remark}
For comparison, we give Gr\"ochenig's notion of localization from
\cite{Gro04}.
Let $I$ and $J$ be countable index sets in $\R^d$ that are separated,
i.e., $\inf_{i \ne j \in I} |i-j| > 0$ and similarly for $J$.
Then $\Fc = \set{f_i}_{i \in I}$ is \emph{$s$-polynomially localized}
with respect to a Riesz basis $\Ec = \set{e_j}_{j \in J}$ if
for every $i \in I$ and $j \in J$ we have
$$|\ip{f_i}{e_j}| \LE C \, (1 + |i-j|)^{-s}
\qquad\text{and}\qquad
|\ip{f_i}{\te_j}| \LE C \, (1 + |i-j|)^{-s},$$
where $\set{\te_j}_{j \in J}$ is the dual basis to 
$\set{e_j}_{j \in J}$.
Likewise $\Fc = \set{f_i}_{i \in I}$ is \emph{exponentially localized}
with respect to a Riesz basis $\Ec = \set{e_j}_{j \in J}$ if
for some $\alpha > 0$ we have for every $i \in I$ and $j \in J$ that
$$|\ip{f_i}{e_j}| \LE C \, e^{-\alpha |i-j|}
\qquad\text{and}\qquad
|\ip{f_i}{\te_j}| \LE C \, e^{-\alpha |i-j|}.
\quad\qed$$
\end{remark}

\subsection{The Approximation Properties} \label{approxsec}

In this section we introduce a collection of definitions which we call
approximation properties.
These definitions extract the essence of the Homogeneous Approximation
Property that is satisfied by Gabor frames, but without reference to the
exact structure of Gabor frames.
A weak HAP for Gabor frames was introduced in \cite{RS95} and developed
further in \cite{GR96}, \cite{CDH99}.
In those papers, the HAP was stated in a form that is specific to the
particular structure of Gabor frames or windowed exponentials,
whereas the following definition applies to arbitrary frames.

\begin{definition}[Homogeneous Approximation Properties] \label{approxdef}
Let $\Fc = \set{f_i}_{i \in I}$ be a frame for~$H$ with
canonical dual $\tFc = \set{\tf_i}_{i \in I}$, and
let $\Ec = \set{e_j}_{j \in G}$ be a sequence in $H$.
Let $a \colon I \to G$ be an associated map.

\smallskip
\begin{enumerate}
\item[(a)]
We say $(\Fc,a,\Ec)$ has the \emph{weak HAP}
if for every $\eps > 0$, there is an integer $N_\eps > 0$
so that for every $j \in G$ we have
$$\dist\Bigparen{e_j, \,\,
\clspan\bigset{\tf_i : i \in I_{N_\eps}(j)}}
\LT \eps.$$
Equivalently, there must exist scalars $c_{i,j}$, with only finitely many
nonzero, such that
\begin{equation} \label{weakHAPdef}
\Bignorm{e_j - \sum_{i \in I_{N_\eps}(j)} c_{i,j} \, \tf_i}
\LT \eps.
\end{equation}

\medskip
\item[(b)]
We say $(\Fc,a,\Ec)$ has the \emph{strong HAP}
if for every $\eps > 0$, there is an integer $N_\eps > 0$
so that for every $j \in G$ we have
\begin{equation} \label{strongHAPdef}
\Bignorm{e_j - \sum_{i \in I_{N_\eps}(j)} \ip{e_j}{f_i} \, \tf_i}
\LT \eps.
\quad\qed
\end{equation}
\end{enumerate}
\end{definition}

We could also define the weak and strong HAPs for frame sequences.
If $\Fc$ is a frame sequence, then a necessary condition for
\eqref{weakHAPdef} or \eqref{strongHAPdef} to hold is that
$\clspan(\Ec) \subset \clspan(\Fc)$.
Thus, the HAPs for frame sequences are the same as the HAPs for a frame
if we set $H = \clspan(\Fc)$.

We also introduce the following symmetric version of the HAPs.

\begin{definition}[Dual Homogeneous Approximation Properties]
Let $\Fc = \set{f_i}_{i \in I}$ be a sequence in $H$, and
let $\Ec = \set{e_j}_{j \in G}$ be a frame for $H$ with
canonical dual $\tEc = \set{\te_j}_{j \in G}$.
Let $a \colon I \to G$ be an associated map.

\smallskip
\begin{enumerate}
\item[(a)]
We say $(\Fc,a,\Ec)$ has the \emph{weak dual HAP}
if for every $\eps > 0$, there is an integer $N_\eps > 0$
so that for every $i \in I$ we have
$\dist\bigparen{f_i, \,\,
\clspan\bigset{\te_j : j \in S_{N_\eps}(a(i))}}
< \eps$.

\medskip
\item[(b)]
We say $(\Fc,a,\Ec)$ has the \emph{strong dual HAP}
if for every $\eps > 0$, there is an integer $N_\eps > 0$
so that for every $i \in I$ we have
$\bignorm{f_i - \sum_{j \in S_{N_\eps}(a(i))} \ip{f_i}{e_j} \, \te_j}
< \eps$.
~\qed
\end{enumerate}
\end{definition}

\subsection{Relations Among the Localization and Approximation Properties}
\label{relationsec}

The following theorem summarizes the relationships that hold among the
localization and approximation properties.
This result is proved in Part~II.

\begin{theorem} \label{relations}
Let $\Fc = \set{f_i}_{i \in I}$ and
$\Ec = \set{e_j}_{j \in G}$ be sequences in $H$, and
let $a : I \to G$ be an associated map.
Then the following statements hold.

\smallskip
\begin{enumerate}
\item[(a)] If $\Fc$ is a frame for $H$, then
$\ell^2$-column decay implies the strong HAP.

\smallskip
\item[(b)] If $\Fc$ is a frame for $H$ and
$\sup_j \norm{e_j} < \infty$, then
the strong HAP implies $\ell^2$-column decay.

\smallskip
\item[(c)] If $\Ec$ is a frame for $H$, then
$\ell^2$-row decay implies the strong dual HAP.

\smallskip
\item[(d)] If $\Ec$ is a frame for $H$ and
$\sup_i \norm{f_i} < \infty$, then
the strong dual HAP implies $\ell^2$-row decay.

\smallskip
\item[(e)] If $\Fc$ is a frame for $H$, then
the strong HAP implies the weak HAP.
If $\Fc$ is a Riesz basis for $H$, then
the weak HAP implies the strong HAP.

\smallskip
\item[(f)] If $\Ec$ is a frame for $H$, then
the strong dual HAP implies the weak dual HAP.
If $\Ec$ is a Riesz basis for $H$, then
the weak dual HAP implies the strong dual HAP.

\smallskip
\item[(g)] If $D^+(I,a) < \infty$ and $1 \le p < \infty$, then
$\ell^p$-localization implies both $\ell^p$-column and $\ell^p$-row decay.
\end{enumerate}
\end{theorem}

\begin{figure}[ht]
\scalebox{.6}{\includegraphics{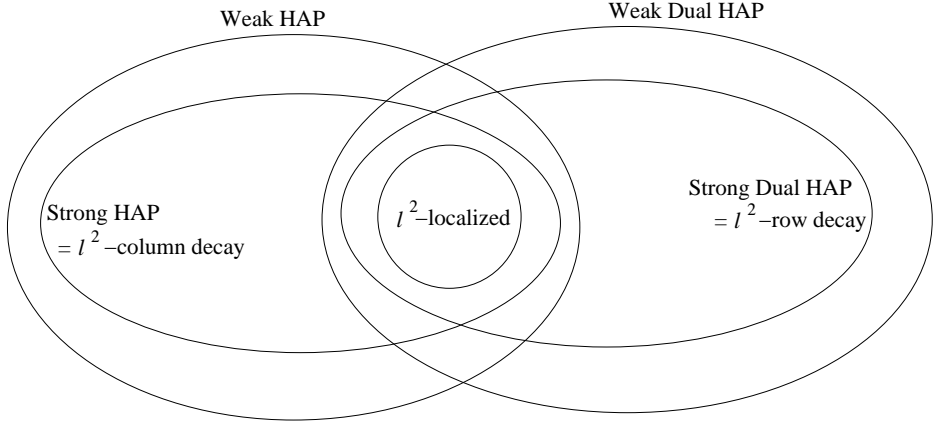}}
\caption{Relations among the localization and approximation properties
for $p=2$, under the assumptions that $\Fc$, $\Ec$ are frames and
$D^+(I,a)<\infty$. \label{fig1}}
\end{figure}

For the case that $\Fc$ and $\Ec$ are both frames for $H$ and 
the upper density $D^+(I,a)$ is finite,
these relations can be summarized in the diagram in Figure~\ref{fig1}.

Part~II exhibits counterexamples to most of
the converse implications of Theorem~\ref{relations}.
These are summarized below.

\begin{enumerate}
\item[(a)]
There exist orthonormal bases $\Ec$, $\Fc$ such that
$(\Fc,a,\Ec)$ does not have $\ell^2$-column decay, and hence does
not satisfy the strong HAP.

\smallskip
\item[(b)]
There exists a frame $\Fc$ and orthonormal basis $\Ec$ such that
$(\Fc,a,\Ec)$ satisfies the weak HAP but not the strong HAP.

\smallskip
\item[(c)]
There exists a frame $\Fc$ and orthonormal basis $\Ec$ such that
$D^+(I,a) < \infty$, $(\Fc,a,\Ec)$ has both $\ell^2$-column decay and
$\ell^2$-row decay, but fails to have $\ell^2$-localization.

\smallskip
\item[(d)]
There exists a Riesz basis $\Fc$ and orthonormal basis $\Ec$ such that
$(\Fc,a,\Ec)$ has $\ell^2$-column decay but not $\ell^2$-row decay.
\end{enumerate}

\subsection{Self-Localization} \label{selflocsec}

In this section we introduce a type of localization in which the system
$\Fc = \set{f_i}_{i \in I}$ is compared to itself or to its canonical
dual frame instead of to a reference system $\Ec$.
An analogous polynomial or exponential ``intrinsic localization'' was
independently introduced by Gr\"ochenig in \cite{Gro03};
see also \cite{For03}, \cite{GF04}.
Although there is no reference system, we still require a mapping
$a \colon I \to G$ associating $I$ with a group $G$.

\begin{definition}[Self-localization]
Let $\Fc = \set{f_i}_{i \in I}$ be a sequence in $H$, and
let $a \colon I \to G$ be an associated map.

\smallskip
\begin{enumerate}
\item[(a)]
We say that $(\Fc,a)$ is \emph{$\ell^p$-self-localized} if there exists
$r \in \ell^p(G)$ such that
$$\forall\, i, j \in I, \quad
|\ip{f_i}{f_j}| \LE r_{a(i) - a(j)}.$$

\medskip
\item[(b)]
If $\Fc$ is a frame sequence, then we say that $(\Fc,a)$ is
\emph{$\ell^p$-localized with respect to its canonical dual frame sequence}
$\tFc = \set{\tf_i}_{i \in I}$ if there exists $r \in \ell^p(G)$ such that
$$\forall\, i, j \in I, \quad
|\ip{f_i}{\tf_j}| \LE r_{a(i) - a(j)}. \quad\qed$$
\end{enumerate}
\end{definition}

\begin{remark} \label{selfremark}
(a) If $I = G$ and $a$ is the identity map, then $(\Fc,a)$ is
$\ell^1$-self-localized if and only if $(\Fc,a,\Fc)$ is $\ell^1$-localized.
However, if $a$ is not the identity map, then this need not be the case.
For example, every orthonormal basis is $\ell^1$-self-localized regardless
of which map~$a$ is chosen, but in Part~II we construct
an orthonormal basis $\Fc = \set{f_i}_{i \in \Z}$
and a map $a \colon \Z \to \Z$ such that $(\Fc,a,\Ec)$ is not
$\ell^1$-localized for any Riesz basis $\Ec$; in fact,
$(\Fc,a,\Ec)$ cannot even possess both $\ell^2$-column decay and
$\ell^2$-row decay for any Riesz basis $\Ec$.

\medskip
(b) Let $\Fc$ be a frame, $\tFc$ its canonical dual frame,
and $S^{-1/2}(\Fc)$ its canonical Parseval frame.
Since $\ip{f_i}{\tf_j} = \ip{S^{-1/2}f_i}{S^{-1/2}f_j}$,
we have that
$(\Fc,a)$ is $\ell^p$-localized with respect to its canonical dual frame
if and only if $(S^{-1/2}(\Fc),a)$ is $\ell^p$-self-localized.
\qed
\end{remark}

We show in Part~II that $\ell^1$-localization with 
respect to the dual frame does not imply $\ell^1$-self-localization.
However, the following result states that the converse is true.
The proof of this result requires us to develop some results on the
Banach algebra of matrices with $\ell^1$-type decay, and is presented in
Appendix~\ref{selflocappend}.
In particular, the proof requires an application of
a type of noncommutative Wiener's Lemma (Theorem~\ref{sjostrand}).

\begin{theorem} \label{selflocthm}
Let $\Fc = \set{f_i}_{i \in I}$ be a frame for $H$, and let
$a \colon I \to G$ be an associated map such that $D^+(I,a) < \infty$.
Let $\tFc$ be the canonical dual frame and
$S^{-1/2}(\Fc)$ the canonical Parseval frame.
If $(\Fc,a)$ is $\ell^1$-self-localized, then:

\begin{enumerate}
\item[(a)]
$(\Fc,a)$ is $\ell^1$-localized with respect to its
canonical dual frame $\tFc = \set{\tf_i}_{i \in I}$,

\medskip
\item[(b)]
$(\tFc,a)$ is $\ell^1$-self-localized, and

\medskip
\item[(c)]
($S^{-1/2}(\Fc),a)$ is $\ell^1$-self-localized.
\end{enumerate}
\end{theorem}

The following is a useful lemma on the relation between self-localization
and localization with respect to a reference sequence.

\begin{lemma} \label{selflemma}
Let $\Fc = \set{f_i}_{i \in I}$ be a sequence in $H$.
Let $\Ec = \set{e_j}_{j \in G}$ be a frame for~$H$ with canonical
dual frame $\tEc$.
Let $a : I \to G$ be an associated map.
If $(\Fc,a,\Ec)$ and $(\Fc,a,\tEc)$ are both $\ell^1$-localized, then
$(\Fc,a)$ is $\ell^1$-self-localized.
In particular, if $\Ec$ is a tight frame and $(\Fc,a,\Ec)$ is
$\ell^1$-localized, then $(\Fc,a)$ is $\ell^1$-self-localized.
\end{lemma}

\begin{proof}
By definition, there exists $r \in \ell^1(G)$ such that both
$|\ip{f_i}{e_j}| \le r_{a(i) - j}$
and
$|\ip{f_i}{\te_j}| \le r_{a(i) - j}$
hold for all $i \in I$ and $j \in G$.
Let $\tilde r(k) = r(-k)$. Then
$$|\ip{f_i}{f_j}|
\EQ \biggabs{\sum_{k \in G} \ip{f_i}{e_k} \, \ip{\te_k}{f_j}}
\LE \sum_{k \in G} r_{a(i) - k}  \, r_{a(j) - k}
\EQ (r * \tilde r)_{a(i) - a(j)}.$$
Since $r * \tilde r \in \ell^1(G)$, we conclude that $(\Fc,a)$ is
$\ell^1$-self-localized.
\end{proof}

\subsection{Relative Measure} \label{relativemeasure}

We now define the relative measure of frame sequences.

\begin{definition} \label{redundancydef}
Let $\Fc = \set{f_i}_{i \in I}$ and $\Ec = \set{e_j}_{j \in G}$
be frame sequences in $H$, and let 
$a \colon I \to G$ be an associated map.
Let $P_\Fc$, $P_\Ec$ denote the orthogonal projections of $H$ onto
$\clspan(\Fc)$ and $\clspan(\Ec)$, respectively.
Then given a free ultrafilter $p$ and a sequence of centers
$c = (c_N)_{N \in \N}$ in $G$, we define the
\emph{relative measure of $\Fc$ with respect to $\Ec$, $p$, and $c$} to be
$$\cM_\Ec(\Fc; p, c)
\EQ \plim_{N \in \N} \frac1{|I_N(c_N)|} \sum_{i \in I_N(c_N)}
    \ip{P_\Ec f_i}{\tf_i}.$$
The \emph{relative measure of $\Ec$ with respect to $\Fc$} is
$$\cM_\Fc(\Ec; p, c)
\EQ \plim_{N \in \N} \frac1{|S_N(c_N)|} \sum_{j \in S_N(c_N)}
      \ip{P_\Fc \te_j}{e_j}.
\quad\qed$$
\end{definition}

Let $A$, $B$ be frame bounds for $\Fc$, and let $E$, $F$ be frame bounds
for $\Ec$.
Then we have the estimates
$|\ip{P_\Ec f_i}{\tf_i}| \le \norm{f_i} \, \norm{\tf_i} \le \sqrt{B/A}$ and
$|\ip{P_\Fc e_j}{\te_j}| \le \norm{e_j} \, \norm{\te_j} \le \sqrt{F/E}$.
Thus,
$|\cM_\Ec(\Fc; p, c)| \le \sqrt{B/A}$ and
$|\cM_\Fc(\Ec; p, c)| \le \sqrt{F/A}$.
Unfortunately, in general $\cM_\Ec(\Fc; p, c)$ or
$\cM_\Fc(\Ec; p, c)$ need be real.
However, if the closed span of $\Fc$ is included in the closed span of $\Ec$
then, as noted in the following definition, the relative measure of $\Ec$
with respect to $\Fc$ will be real and furthermore we can give tighter
bounds on its value, as pointed out in the following definition.

\begin{definition}
If $\clspan(\Ec) \supset \clspan(\Fc)$ then $P_\Ec$ is the identity map and
$\Ec$ plays no role in determining the value of $\cM_\Ec(\Fc;p,e)$.
Therefore, in this case we define the
\emph{measure of $\Fc$ with respect to $p$ and $c$} to be
$$\cM(\Fc; p, c)
\EQ \plim_{N \in \N} \frac1{|I_N(c_N)|} \sum_{i \in I_N(c_N)}
    \ip{f_i}{\tf_i}.$$
Since $\ip{f_i}{\tf_i} = \norm{S^{-1/2}f_i}^2$, we have that
$\cM(\Fc; p, c)$ is real.
Additionally, since $S^{-1/2}(\Fc)$ is a Parseval frame, we have
$0 \le \ip{f_i}{\tf_i} \le 1$ for all $i$, and therefore
$$0 \LE \cM(\Fc;p,c) \LE 1.$$
We further define the 
\emph{lower and upper measures of $\Fc$} to be, respectively,
\begin{align}
\cM^-(\Fc)
& \EQ \liminf_{N \to \infty} \, \inf_{j \in G} \,
      \frac1{|I_N(j)|} \sum_{i \in I_N(j)}
      \ip{f_i}{\tf_i}, \label{lowermeasure}
      \\[1 \jot]
\cM^+(\Fc)
& \EQ \limsup_{N \to \infty} \, \sup_{j \in G} \,
      \frac1{|I_N(j)|} \sum_{i \in I_N(j)}
      \ip{f_i}{\tf_i}. \label{uppermeasure}
\end{align}
As in Lemma~\ref{liminfexist}, there will exist free ultrafilters
$p^-$, $p^+$ and sequence of centers $c^-$, $c^+$ such that
$\cM^-(\Fc) = \cM(\Fc;p^-,c^-)$ and
$\cM^+(\Fc) = \cM(\Ec;p^+,c^+)$.

When $\clspan(\Fc) \supset \clspan(\Ec)$, we define
$\cM(\Ec; p, c)$ and $\cM^\pm(\Ec)$ in an analogous manner.
~\qed
\end{definition}

\begin{example} \label{specialcases}
The following special cases show that the measure of a Riesz basis is~$1$.

\smallskip
\begin{enumerate}
\item[(a)] If $\clspan(\Ec) \supset \clspan(\Fc)$ and
$\Fc$ is a Riesz sequence then
$\ip{f_i}{\tf_i} = 1$ for every $i \in I$, so
$\cM(\Fc;p,c) = \cM^+(\Fc) = \cM^-(\Fc) = 1$.

\medskip
\item[(b)] If $\clspan(\Fc) \supset \clspan(\Ec)$ and
$\Ec$ is a Riesz sequence then
$\ip{\te_j}{e_j} = 1$ for every $j \in G$, so
$\cM(\Ec;p,c) = \cM^+(\Ec) = \cM^-(\Ec) = 1$.
~\qed
\end{enumerate}
\end{example}

\begin{example}
For each $k = 1, \dots, M$, let $\set{f_{jk}}_{j \in \Z}$ be an
orthogonal basis for~$H$ such that
$\norm{f_{jk}}^2 = A_k$ for every $j \in \Z$.
Let $I = \Z \times \set{1,\dots,M}$.
Then $\Fc = \set{f_{jk}}_{(j,k) \in I}$ is a tight frame for $H$
and its canonical dual frame is
$\tFc = \set{\tf_{jk}}_{(j,k) \in I}$
where $\tf_{jk} = (\frac1{A_1 + \cdots + A_M}) \, f_{jk}$.
Define $a \colon I \to \Z$ by $a(j,k) = j$.
Then for each $N$,
$$\frac1{|I_N(c_N)|} \sum_{(j,k) \in I_N(c_N)} \ip{f_{jk}}{\tf_{jk}}
\EQ \frac1{MN} \sum_{k=1}^M \, \sum_{j \in [c_N - \frac{N}2, c_N + \frac{N}2)}
    \frac{A_k}{A_1 + \cdots + A_M}
\EQ \frac1M.$$
Consequently, for any choice of free ultrafilter $p$ or sequence of centers $c$
we have
$\cM(\Fc;p,c) = \cM^-(\Fc) = \cM^+(\Fc) = \frac1M$.
~\qed
\end{example}

\begin{example}[Lattice Gabor Systems] \label{latticegabor}
Consider a lattice Gabor frame, i.e., a frame of the form
$\Gc(g,\alpha\Z^d \times \beta\Z^d)$.
The canonical dual frame is a lattice Gabor frame of the form
$\Gc(\tg,\alpha\Z^d \times \beta\Z^d)$ for some $\tg \in L^2(\R^d)$.
By the Wexler--Raz relations, we have $\ip{g}{\tg} = (\alpha\beta)^d$
(we also derive this fact directly from our results in Part~II).
Since
$\ip{M_{\beta n} T_{\alpha k} g}{M_{\beta n} T_{\alpha k} \tg}
= \ip{g}{\tg}$,
we therefore have for any free ultrafilter $p$ and sequence of centers
$c = (c_N)_{N \in \N}$ in $\alpha\Z^d \times \beta\Z^d$ that
$$\cM(\Gc(g,\alpha\Z^d \times \beta\Z^d);p,c)
\EQ \cM^\pm(\Gc(g,\alpha\Z^d \times \beta\Z^d))
\EQ \ip{g}{\tg}
\EQ (\alpha\beta)^d.$$
Since we also have
$\BD^\pm(\alpha\Z^d \times \beta\Z^d) = (\alpha\beta)^{-d}$,
we conclude that
$$\cM^\pm(\Gc(g,\alpha\Z^d \times \beta\Z^d))
\EQ \frac1{\BD^\mp(\alpha\Z^d \times \beta\Z^d)}.$$
We prove a similar but much more general relationship for
abstract localized frames in Theorems~\ref{densityredundancy}
and~\ref{excesscor}.
~\qed
\end{example}

The following proposition gives a connection between measure and
excess (excess was defined just prior to equation \eqref{excesscondition}).
By imposing localization hypotheses, stronger results will be derived
in Section~\ref{positivesec}.

\begin{proposition}[Infinite Excess] \label{infiniteexcess}
Let $\Fc = \set{f_i}_{i \in I}$ be a frame sequence and $a \colon I \to G$
an associated map.
If $\cM^-(\Fc) < 1$, then $\Fc$ has infinite excess, and furthermore,
there exists an infinite subset $J \subset I$ such that
$\set{f_i}_{i \in I \setminus J}$ is still a frame for $\clspan(\Fc)$.
\end{proposition}
\begin{proof}
Fix $s$ with $\cM^-(\Fc) < s < 1$.
Then, considering the definition of $\cM^-(\Fc)$ in \eqref{lowermeasure},
there exists a subsequence $N_k \to \infty$ and points $j_k$ such that
$$\frac1{|I_{N_k}(j_k)|} \sum_{i \in I_{N_k}(j_k)} \ip{f_i}{\tf_i}
\LE s
\LT 1$$
for each $k$.
It then follows that there exists an infinite subset
$J \subset I$ such that $\sup_{i \in J} \ip{f_i}{\tf_i} < 1$,
which by \cite[Cor.~5.7]{BCHL03} completes the proof.
\end{proof}

In general, the set $J$ constructed in the preceding proposition may have
zero density.
The following result provides a necessary condition under which a set of
positive density can be removed yet leave a frame
(a sufficient condition will be obtained in
Theorem~\ref{positiveremove} below).
For simplicity of notation, if $J \subset I$ then we will write
$D(p,c;J,a)$ to mean $D(p,c;J,a|_J)$.

\begin{proposition} \label{jalphaprop}
Let $\Fc = \set{f_i}_{i \in I}$ be a frame sequence and $a \colon I \to G$
an associated map such that $0 < D^-(I,a) \le D^+(I,a) < \infty$.
For each $0 \le \alpha \le 1$, define
\begin{equation} \label{Jalphadef}
J_\alpha
\EQ \set{i \in I : \ip{f_i}{\tf_i} \le \alpha}.
\end{equation}
Then the following statements hold.

\smallskip
\begin{enumerate}
\item[(a)]
For each free ultrafilter $p$ and sequence of centers $c = (c_N)_{N \in \N}$
in $G$, we have for each $0 < \alpha < 1$ that
\begin{align}
\frac{\alpha - \cM(\Fc;p,c)}{\alpha} \, D(p,c;I,a)
& \LE D(p,c;J_{\alpha},a) \label{Jalphaest1} \\
& \LE \frac{1- \cM(\Fc;p,c)}{1 - \alpha} \, D(p,c;I,a). \label{Jalphaest2}
\end{align}

\medskip
\item[(b)]
If there exists a free ultrafilter $p$ and sequence of centers
$c = (c_N)_{N \in \N}$ in $G$ such that
$D(p,c;J_\alpha,a) > 0$, then $\cM(\Fc;p,c) < 1$.
Consequently $\cM^-(\Fc) < 1$ and there exists an infinite set $J \subset I$
such that $\set{f_i}_{i \in I \setminus J}$ is a frame for $\clspan(\Fc)$.

\medskip
\item[(c)]
If there exists a
subset $J \subset I$, a free ultrafilter $p$, and a sequence of centers
$c = (c_N)_{N \in \N}$ in $G$ such that $D(p,c;J,a) > 0$ and
$\set{f_i}_{i \in I \setminus J}$ is a frame for $\clspan(\Fc)$,
then $\cM(\Fc;p,c) < 1$.
In particular, $\cM^-(\Fc) < 1$.
\end{enumerate}
\end{proposition}

\begin{proof}
(a) Consider any $0 < \alpha < 1$.
If $\cM(\Fc;p,c) \ge \alpha$ then inequality \eqref{Jalphaest1}
is trivially satisfied, so assume that $\cM(\Fc;p,c) < \alpha$.
Fix $\eps > 0$ so that $\cM(\Fc;p,c) + \eps < \alpha$.
Then by definition of ultrafilter, there exists an infinite set $A \in p$
such that
\begin{equation} \label{limD}
\forall\, N \in A, \quad
\Bigabs{\cM(\Fc;p,c) - \frac1{|I_N(c_N)|}
\sum_{i \in I_N(c_N)} \ip{f_i}{\tf_i}}
\LT \eps.
\end{equation}
Hence for $N \in A$ we have
\begin{align*}
\cM(\Fc;p,c) + \eps
& \GE \frac1{|I_N(c_N)|} \sum_{i \in I_N(c_N)} \ip{f_i}{\tf_i} \\[1 \jot]
& \EQ \frac1{|I_N(c_N)|} \,
      \biggparen{\sum_{i \in I_N(c_N) \cap J_\alpha}  \ip{f_i}{\tf_i} \plus
      \sum_{i \in I_N(c_N) \cap J_\alpha^\comp}  \ip{f_i}{\tf_i}}
      \allowdisplaybreaks \\[1 \jot]
& \GE \frac{0 \cdot |I_N(c_N) \cap J_\alpha| \plus
      \alpha \cdot |I_N(c_N) \cap J_\alpha^\comp|} {|I_N(c_N)|} \\[1 \jot]
& \EQ \alpha \, \frac{|I_N(c_N)| - |I_N(c_N) \cap J_\alpha|} {|I_N(c_N)|}.
\end{align*}
Multiplying both sides of this inequality by $\frac{|I_N(c_N)|}{|S_N(c_N)|}$
and rearranging, we find that
$$\forall\, N \in A, \quad
\frac{|I_N(c_N) \cap J_\alpha|} {|S_N(c_N)|}
\GE \biggparen{1 - \frac{\cM(\Fc;p,c) + \eps}{\alpha}} \,
    \frac{|I_N(c_N)|} {|S_N(c_N)|}.$$
Taking the limit with respect to the ultrafilter $p$ we obtain
$$D(p,c;J_\alpha,a)
\GE \biggparen{1 - \frac{\cM(\Fc;p,c) + \eps}{\alpha}} \, D(p,c;I,a).$$
Since $\eps$ was arbitrary, we obtain the inequality \eqref{Jalphaest1}.

The inequality \eqref{Jalphaest2} is similar,
arguing from an infinite set $A \in p$ such that \eqref{limD} holds true that
\begin{align*}
\cM(\Fc;p,e) - \eps
& \LE \frac1{|I_N(c_N)|} \sum_{i \in I_N(c_N)} \ip{f_i}{\tf_i} \\[1 \jot]
& \EQ \frac1{|I_N(c_N)|} \,
      \biggparen{\sum_{i \in I_N(c_N) \cap J_\alpha}  \ip{f_i}{\tf_i} \plus
      \sum_{i \in I_N(c_N) \cap J_\alpha^\comp}  \ip{f_i}{\tf_i}}
      \allowdisplaybreaks \\[1 \jot]
& \LE \frac{\alpha \cdot |I_N(c_N) \cap J_\alpha| \plus
      1 \cdot |I_N(c_N) \cap J_\alpha^\comp|} {|I_N(c_N)|} \\[1 \jot]
& \EQ \frac{|I_N(c_N)| -
      (1-\alpha) \cdot |I_N(c_N) \cap J_\alpha|} {|I_N(c_N)|},
\end{align*}
and then multiplying both sides of this inequality by
$\frac{|I_N(c_N)|}{|S_N(c_N)|}$,
rearranging, and taking the limit.

\medskip
(b) Follows immediately from~(a) and Proposition~\ref{infiniteexcess}.

\medskip
(c) Suppose that such a $J$ exists.
If $f_i = 0$ for every $i \in J$ then the result is trivial, so suppose this
is not the case.
Let $S$ be the frame operator for~$\Fc$.
Then $\set{S^{-1/2}f_i}_{i \in I \setminus J}$ is a frame, and in
particular is a subset of the Parseval frame $S^{-1/2}(\Fc)$.
For a given $j \in J$, the optimal lower frame bound for the frame
$\set{S^{-1/2}f_i}_{i \ne j}$ with a single element deleted is
$1 - \norm{S^{-1/2}f_j}^2 = 1 - \ip{f_j}{\tf_j}$.
Hence, if $A$ is a lower frame bound for
$\set{S^{-1/2}f_i}_{i \in I \setminus J}$,
then $A \le 1 - \ip{f_j}{\tf_j}$ for all $j \in J$.
Thus $J \subset J_\alpha$ where $\alpha = 1 - A$, and consequently,
for any $p$ and $c$ we have $D(p,c;J_\alpha,a) \ge D(p,c;J,a) > 0$.
Therefore \eqref{Jalphaest2} implies that $\cM(\Fc;p,c) < 1$.
\end{proof}

Choosing in the preceding proposition the ultrafilters $p$ and centers $c$ 
that achieve upper or lower density or measure yields the following corollary.

\begin{corollary} \label{Jalphacoro}
Let $\Fc = \set{f_i}_{i \in I}$ be a frame sequence and $a \colon I \to G$
an associated map such that $0 < D^-(I,a) \le D^+(I,a) < \infty$.
Let $J_\alpha$ be defined by \eqref{Jalphadef}.
Then the following statements hold.

\smallskip
\begin{enumerate}
\item[(a)]
$\cM^+(\Fc) < 1$ if and only if there exists $0 < \alpha < 1$
such that $D^-(J_\alpha,a) > 0$.
In fact, $D^-(J_\alpha,a) > 0$ for all
$\cM^+(\Fc) < \alpha < 1$.

\medskip
\item[(b)]
If there exists $J \subset I$ such that $D^-(J,a) > 0$ and
$\set{f_i}_{i \in I \setminus J}$ is a frame for $\clspan(\Fc)$,
then $\cM^+(\Fc) < 1$.

\medskip
\item[(c)]
$\cM^-(\Fc) < 1$ if and only if there exists $0 < \alpha < 1$
such that $D^+(J_\alpha,a) > 0$.
In fact, $D^+(J_\alpha,a) > 0$ for all
$\cM^-(\Fc) < \alpha < 1$.

\medskip
\item[(d)]
If there exists $J \subset I$ such that $D^+(J,a) > 0$ and
$\set{f_i}_{i \in I \setminus J}$ is a frame for $\clspan(\Fc)$,
then $\cM^-(\Fc) < 1$.
\end{enumerate}
\end{corollary}
\begin{proof}
Suppose that $\cM^+(\Fc) < 1$, and fix
$\cM^+(\Fc) < \alpha < 1$.
Let $p$ and $c$ be the free ultrafilter and sequence of centers given
by Lemma~\ref{liminfexist}(b) such that
$D^-(J_\alpha,a) = D(p,c;J_\alpha,a)$.
Then by Proposition~\ref{jalphaprop},
\begin{align*}
D^-(J_\alpha,a)
\EQ D(p,c;J_\alpha,a)
& \GE \frac{\alpha - \cM(F;p,c)}{\alpha} \, D(p,c;I,a) \\[1 \jot]
& \GE \frac{\alpha - \cM^+(F)}{\alpha} \, D^-(I,a)
\GT 0.
\end{align*}
The other statements are similar.
\end{proof}

\section{Density and Overcompleteness} \label{mainsec}

\subsection{Necessary Density Conditions} \label{necessarysec}

In this section we prove two necessary conditions on the density of
localized frames.

First we require the following standard lemma.

\begin{lemma} \label{dimcount}
Let $H_N$ be an $N$-dimensional Hilbert space.
Then the following statements hold.

\begin{enumerate}
\item[(a)]
Let nonzero $f_1,\dots,f_M \in H_N$ be given.
Let $m = \min\set{\norm{f_1},\dots,\norm{f_M}}$.
Then the Bessel bound $B$ for $\set{f_1,\dots,f_M}$ satisfies
$B \ge mM/N$.

\medskip
\item[(b)] If $\set{f_i}_{i \in J}$ is a Bessel sequence in $H_N$
that is norm-bounded below, i.e., $\inf_i \norm{f_i} > 0$, then $J$ is finite.
\end{enumerate}
\end{lemma}
\begin{proof}
(a) We may assume that $H_N = \Span\set{f_1,\dots,f_M}$.
Then $\set{f_1,\dots,f_M}$ is a frame for $H_N$, so this family has a positive
definite frame operator $S$.
Let $\lambda_1 \ge \dots \ge \lambda_N$ be the eigenvalues of $S$.
Letting $\set{\tf_1,\dots,\tf_M}$ be the dual frame, we have then that
$$\sum_{j=1}^N \lambda_j
\EQ \trace(S)
\EQ \sum_{i=1}^M \ip{Sf_i}{\tf_i}
\EQ \sum_{i=1}^M \norm{f_i}^2
\GE mM.$$
Hence $mM/N \le \lambda_1 = \norm{S} \le B$.

\medskip
(b) From part~(a), $|J| \le BN/m < \infty$.
\end{proof}

Our first main result shows that the weak HAP implies a lower bound for the
density of a frame.
The proof is inspired by the double projection techniques of \cite{RS95},
although those results relied on the structure of Gabor frames and, in
particular, a version of the HAP that is satisfied by Gabor frames.

\begin{theorem}[Necessary Density Bounds] \label{necessary} \

\begin{enumerate}
\item[(a)]
Assume $\Fc = \set{f_i}_{i \in I}$ is a frame for $H$ 
and $\Ec = \set{e_j}_{j \in G}$ is a Riesz sequence in $H$.
Let $a \colon I \to G$ be an associated map.
If $(\Fc,a,\Ec)$ has the weak HAP, then 
$$1 \LE D^-(I,a) \LE D^+(I,a) \LE \infty.$$

\medskip
\item[(b)]
Assume $\Fc = \set{f_i}_{i \in I}$ is a Riesz sequence in $H$ 
and $\Ec = \set{e_j}_{j \in G}$ is a frame for $H$.
Let $a \colon I \to G$ be an associated map.
If $(\Fc,a,\Ec)$ has the weak dual HAP, then 
$$0 \LE D^-(I,a) \LE D^+(I,a) \LE 1.$$
\end{enumerate}
\end{theorem}
\begin{proof}
(a) Let $\tFc = \set{\tf_i}_{i \in I}$ be the canonical dual frame to $\Fc$,
and let $\tEc = \set{\te_j}_{j \in G}$ be the Riesz sequence in $\clspan(\Ec)$
that is biorthogonal to $\Ec$.
Fix $\eps > 0$, and let $N_\eps$ be the number given in the definition
of the weak HAP.
Fix an arbitrary point $j_0 \in G$ and a box size $N > 0$.
Define
$$V \EQ \Span\set{e_j : j \in S_N(j_0)}
\qquad\text{and}\qquad
W \EQ \clspan\set{\tf_i : i \in I_{N+N_\eps}(j_0)}.$$
Note that $V$ is finite-dimensional, with $\dim(V) = |S_N(j_0)|$.
On the other hand, $W$ may be finite or infinite-dimensional, but in any case
we have
$\Dim(W) \le |I_{N+N_\eps}(j_0)|$ in the sense of the extended reals.

Let $P_V$ and $P_W$ denote the orthogonal projections of $H$ onto $V$ and $W$,
respectively.
Define a map $T \colon V \to V$ by $T = P_V P_W$.
Note that since the domain of $T$ is $V$, we have
$T = P_V P_W P_V$, so $T$ is self-adjoint.

Let us estimate the trace of $T$.
First note that every eigenvalue $\lambda$ of $T$ satisfies
$|\lambda| \le \norm{T} \le \norm{P_V} \, \norm{P_W} = 1$.
This provides us with an upper bound for the trace of $T$, since
the trace is the sum of the eigenvalues, and hence
\begin{equation} \label{upper}
\trace(T)
\LE \rank(T)
\LE \Dim(W)
\LE |I_{N+N_\eps}(j_0)|.
\end{equation}

For a lower estimate, note that $\set{e_j : j \in S_N(j_0)}$
is a Riesz basis for $V$.
The dual Riesz basis in $V$ is $\set{P_V \te_j : j \in S_N(j_0)}$.
Therefore
\begin{align}
\trace(T)
& \EQ \sum_{j \in S_N(j_0)} \ip{T e_j}{P_V \te_j}
      \label{est1} \\[1 \jot]
& \EQ \sum_{j \in S_N(j_0)} \ip{P_V T e_j}{\te_j}
      \allowdisplaybreaks \notag \\[1 \jot]
& \EQ \sum_{j \in S_N(j_0)} \ip{e_j}{\te_j} \plus
      \sum_{j \in S_N(j_0)} \ip{(P_V P_W - \one) e_j}{\te_j}
      \notag \\[1 \jot]
& \GE |S_N(j_0)| \minus
      \sum_{j \in S_N(j_0)} |\ip{(P_V P_W - \one) e_j}{\te_j}|, \notag
\end{align}
where in the last line we have used the fact that
$\ip{e_j}{\te_j} = 1$.

The elements of any Riesz sequence are uniformly bounded in norm,
so $C = \sup_j \norm{\te_j} < \infty$.
Hence
\begin{equation} \label{est2}
|\ip{(P_V P_W - \one) e_j}{\te_j}|
\LE \norm{(P_V P_W - \one) e_j} \, \norm{\te_j}
\LE C \, \norm{(P_V P_W - \one) e_j}.
\end{equation}
Since $(P_V P_W - \one) e_j \in V$
while $(\one - P_V) P_W e_j \perp V$,
we have by the Pythagorean Theorem that
\begin{align*}
\norm{(P_W - \one) e_j}^2
& \EQ \norm{(P_V P_W - \one) e_j + (\one - P_V) P_W e_j}^2 \\[1 \jot]
& \EQ \norm{(P_V P_W - \one) e_j}^2 + \norm{(\one - P_V) P_W e_j}^2.
\end{align*}
Thus,
\begin{align}
\norm{(P_V P_W - \one) e_j}^2
& \EQ \norm{(P_W - \one) e_j}^2 - \norm{(\one - P_V) P_W e_j}^2
      \label{est3} \\[2 \jot]
& \LE \norm{(P_W - \one) e_j}^2
      \notag \\[2 \jot]
& \EQ \dist(e_j, W)^2. \notag
\end{align}
However, for $j \in S_N(j_0)$, we have
$I_{N_\eps(j)} \subset I_{N+N_\eps}(j_0)$, so for such $j$,
\begin{align}
\dist(e_j, W)
& \EQ \dist\bigparen{e_j, \, \clspan\set{\tf_i : i \in I_{N+N_\eps}(j_0)}}
      \label{est4} \\[1 \jot]
& \LE \dist\bigparen{e_j, \, \clspan\set{\tf_i : i \in I_{N_\eps}(j)}}
\LT \eps, \notag
\end{align}
the last inequality following from the weak HAP.
By combining equations \eqref{est1}--\eqref{est4}, we find that
\begin{equation} \label{lower}
\trace(T)
\GE |S_N(j_0)| - \sum_{j \in S_N(j_0)} C \eps
\EQ (1 - C \eps) \, |S_N(j_0)|.
\end{equation}

Finally, combining the upper estimate for $\trace(T)$ from \eqref{upper}
with the lower estimate from \eqref{lower}, we obtain
$$\frac{|I_{N+N_\eps}(j_0)|} {|S_{N+N_\eps}(j_0)|}
\GE \frac{(1 - C \eps) \, |S_N(j_0)|} {|S_{N+N_\eps}(j_0)|},$$
where the left-hand side could be infinite.
In any case, taking the infimum over all $j_0 \in G$ and then the
liminf as $N \to \infty$ yields
$$D^-(I,a)
\EQ \liminf_{N \to \infty} \inf_{j_0 \in G}
    \frac{|I_{N+N_\eps}(j_0)|}{|S_{N+N_\eps}(j_0)|}
\GE (1 - C \eps) \, \liminf_{N \to \infty}
    \frac{|S_N(j_0)|} {|S_{N+N_\eps}(j_0)|}
\EQ 1 - C \eps,$$
the last equality following from the asymptotics in \eqref{asymptotics}.
Since $\eps$ was arbitrary, we obtain $D^-(I,a) \ge 1$.

\medskip
(b) Let $\tFc = \set{\tf_i}_{i \in I}$ be the Riesz sequence in $\clspan(\Fc)$
that is biorthogonal to $\Fc$, and let
$\tEc = \set{\te_j}_{j \in G}$ be the canonical dual frame to $\Ec$.
Fix $\eps > 0$, and let $N_\eps$ be the number given in the definition
of the weak dual HAP.
Fix an arbitrary point $j_0 \in G$ and a box size $N > 0$.
Define
$$V \EQ \clspan\set{f_i : i \in I_N(j_0)}
\qquad\text{and}\qquad
W \EQ \Span\set{\te_j : j \in S_{N+N_\eps}(j_0)}.$$
Note that $W$ is finite-dimensional, with
$\Dim(W) \le |S_{N+N_\eps}(j_0)|$.
We will show next that $V$ is also finite-dimensional.

Because $\Fc$ is a Riesz sequence, it is norm-bounded below.
In fact, $\norm{f_i} \ge A^{1/2}$ where $A$,~$B$ are frame bounds for $\Fc$.
Now for $i \in I_N(j_0)$ we have
$S_{N_\eps}(a(i)) \subset S_{N+N_\eps}(j_0)$, so
\begin{align*}
\dist(f_i, W)
& \EQ \dist\bigparen{f_i, \, \Span\set{\te_j : j \in S_{N+N_\eps}(j_0)}}
      \\[1 \jot]
& \LE \dist\bigparen{f_i, \, \Span\set{\te_j : j \in S_{N_\eps}(a(i))}}
\LT \eps,
\end{align*}
the last inequality following from the weak dual HAP.
Hence
\begin{equation} \label{lowerbound}
\forall\, i \in I_N(j_0), \quad
\norm{P_W f_i} \GE \norm{f_i} - \eps \GE A^{1/2} - \eps.
\end{equation}
Thus $\set{P_W f_i}_{i \in I_N(j_0)}$
is a Bessel sequence in the finite-dimensional space $W$, and
furthermore this sequence is norm-bounded below by \eqref{lowerbound}.
Lemma~\ref{dimcount} therefore implies that $I_N(j_0)$ is finite.
Thus $V$ is finite-dimensional, as $\dim(V) = |I_N(j_0)| < \infty$.

Let $P_V$ and $P_W$ denote the orthogonal projections of $H$ onto $V$ and $W$,
respectively, and define a map $T \colon V \to V$ by $T = P_V P_W$.
An argument very similar to the one used in part~(a) then shows that
$(1 - C \eps) \, |I_N(j_0)| \le |S_{N+N_\eps}(j_0)|$,
where $C = \sup_i \norm{\tf_i} < \infty$.
Taking the supremum over all $j_0 \in G$ and then the limsup as $N \to \infty$
then yields the result.
\end{proof}

The conclusion of Theorem~\ref{necessary}(a) allows the possibility
that the density might be infinite.
Our next main result will show that $\ell^2$-row decay implies, at least
for Bessel sequences compared to frames, that the upper density is finite.

\begin{theorem}[Necessary Finite Density Condition] \label{finitedensity}
Let $\Fc = \set{f_i}_{i \in I}$ be a Bessel sequence in $H$, and suppose
$\inf_{i \in I} \norm{f_i} > 0$.
Assume $\Ec = \set{e_j}_{j \in G}$ is a frame for~$H$,
and let $a \colon I \to G$ be an associated map.
If $(\Fc,a,\Ec)$ has $\ell^2$-row decay,
then $D^+(I,a) < \infty$.
\end{theorem}
\begin{proof}
If we let $S$ be the frame operator for $\Ec$ then
$S^{-1/2}(\Ec)$ is a Parseval frame for~$H$.
Further, $\ip{f_i}{e_j} = \ip{S^{1/2}f_i}{S^{-1/2}e_j}$ and
$S^{1/2}(\Fc)$ is still a Bessel sequence in $H$ that is norm-bounded below.
Thus, it suffices to show the result when $\Ec$ is a Parseval frame for~$H$.

Let $B$ be the Bessel bound for $\Fc$, and let $m = \inf_i \norm{f_i}^2$.
Fix $0 < \eps < m$.
Since $(\Fc,a,\Ec)$ has $\ell^2$-row decay, there exists an $N_\eps$ such that
$$\forall\, i \in I, \quad
\sum_{j \in G \setminus S_{N_\eps}(a(i))} |\ip{f_i}{e_j}|^2 \LT \eps.$$
Let $j_0 \in G$ and $N > 0$ be given.
Define
$$V \EQ \Span\set{e_j : j \in S_{N+N_\eps}(j_0)},$$
and note that $\dim(V) \leq |S_{N+N_\eps}(j_0)|$.
Define $L_V \colon H \to V$ by
$$L_V f
\EQ \sum_{j \in S_{N+N_\eps}(j_0)} \ip{f}{e_j} \, e_j,
\qquad f \in H,$$
and set $h_i = L_V f_i$ for $i \in I$.
Since $\norm{L_V} \le 1$, it follows that $\set{h_i}_{i \in I}$
is a Bessel sequence in $H$ with the same Bessel bound~$B$ as $\Fc$.

Now, if $i \in I_N(j_0)$ then $a(i) \in S_N(j_0)$, so
$S_{N_\eps}(a(i)) \subset S_{N+N_\eps}(j_0)$.
Therefore,
$$\sum_{j \in G \setminus S_{N+N_\eps}(j_0)} |\ip{f_i}{e_j}|^2
\LE \sum_{j \in G \setminus S_{N_\eps}(a(i))} |\ip{f_i}{e_j}|^2
\LT \eps.$$
Hence
$$\sum_{j \in S_{N+N_\eps}(j_0)} |\ip{f_i}{e_j}|^2
\GE \sum_{j \in G} |\ip{f_i}{e_j}|^2 \minus \eps
\EQ \norm{f_i}^2 - \eps
\GE m - \eps.$$
On the other hand,
$$\sum_{j \in S_{N+N_\eps}(j_0)} |\ip{f_i}{e_j}|^2
\EQ \ip{h_i}{f_i}
\LE \norm{h_i} \, \norm{f_i}
\LE B^{1/2} \, \norm{h_i}.$$
Thus
$$\norm{h_i} \GE \frac{m-\eps}{B^{1/2}},
\qquad i \in I_N(j_0).$$
Applying Lemma~\ref{dimcount}(a) to $\set{h_i}_{i \in I_N(j_0)}$,
we conclude that
$$B \GE \frac{m-\eps}{B^{1/2}} \,\, \frac{|I_N(j_0)|}{\dim(V)}
\GE \frac{m-\eps}{B^{1/2}} \,\, \frac{|I_N(j_0)|}{|S_{N+N_\eps(j_0)}|}.$$
Consequently, applying the asymptotics in \eqref{asymptotics}, we conclude that
\begin{align*}
D^+(I,a)
& \EQ \limsup_{N \to \infty} \, \sup_{j_0 \in G} \,
      \frac{|I_N(j_0)|}{|S_N(j_0)|} \\[1 \jot]
& \LE \limsup_{N \to \infty} \, \sup_{j_0 \in G} \,
      \frac{B^{3/2}}{m-\eps} \, \frac{|S_{N+N_{\eps}}(j_0)|}{|S_N(j_0)|}
\EQ \frac{B^{3/2}}{m-\eps}
\LT \infty. \qedhere
\end{align*}
\end{proof}

\subsection{The Connection Between Density and Relative Measure}
\label{connectsec}

We now derive the fundamental relationship between density
and relative measure for localized frames.

\begin{theorem}[Density--Relative Measure] \label{densityredundancy}
Let $\Fc =\set{f_i}_{i \in I}$ and $\Ec =\set{e_j}_{j \in G}$
be frame sequences in $H$, and let $a \colon I \to G$ be an associated map.
If $D^+(I,a) < \infty$ and $(\Fc,a,\Ec)$
has both $\ell^2$-column decay and $\ell^2$-row decay,
then the following statements hold.

\smallskip
\begin{enumerate}
\item[(a)] For every sequence of centers $c = (c_N)_{N \in \N}$ in $G$,
\begin{align*}
\lim_{N \to \infty} & \biggl[
\biggparen{\frac1{|S_N(c_N)|} \sum_{j \in S_N(c_N)} \ip{P_\Fc \te_j}{e_j}}
\minus \biggr. \\
& \qquad\quad \biggl. \biggparen{\frac{|I_N(c_N)|}{|S_N(c_N)|}} \,
  \biggparen{\frac1{|I_N(c_N)|} \sum_{i \in I_N(c_N)} \ip{P_\Ec f_i}{\tf_i}}
  \biggr]
\EQ 0.
\end{align*}

\smallskip
\item[(b)] For every sequence of centers $c = (c_N)_{N \in \N}$ in $G$ and
any free ultrafilter $p$,
$$\cM_\Fc(\Ec;p,c) \EQ D(p,c) \cdot \cM_\Ec(\Fc;p,c).$$
\end{enumerate}
\end{theorem}

\begin{proof}
(a) Fix any sequence of centers $c = (c_N)_{N \in \N}$ in $G$.
Define
\begin{align*}
d_N & \EQ \frac{|I_N(c_N)|}{|S_N(c_N)|}, \\[2 \jot]
r_N & \EQ \frac1{|I_N(c_N)|} \sum_{i \in I_N(c_N)} \ip{P_\Ec f_i}{\tf_i},
          \\[1 \jot]
s_N & \EQ \frac1{|S_N(c_N)|} \sum_{j \in S_N(c_N)} \ip{P_\Fc \te_j}{e_j}.
\end{align*}
We must show that $|s_N - d_N r_N| \to 0$.

First, we make some preliminary observations and introduce some notation.
Let $A$,~$B$ denote frame bounds for $\Fc$, and let
$E$,~$F$ denote frame bounds for $\Ec$.
Then the canonical dual frame sequences $\tFc$ and $\tEc$ have frame bounds
$\frac1B$,~$\frac1A$ and
$\frac1F$,~$\frac1E$, respectively.
Consequently, for all $i \in I$ and $j \in G$,
$$\norm{f_i}^2 \LE B, \qquad
\norm{\tf_i}^2 \LE \frac1A, \qquad
\norm{e_j}^2 \LE F, \qquad
\norm{\te_j}^2 \LE \frac1E.$$

Fix any $\eps > 0$.
Since $(\Fc,a,\Ec)$ has both $\ell^2$-column decay and $\ell^2$-row decay,
there exists an integer $N_\eps > 0$ such that both
equations \eqref{PropXdef} and \eqref{PropXstardef} hold.
Additionally, since $D^+(I,a) < \infty$, there exists an $K > 0$
such that \eqref{annulus} holds.

Let $P_\Fc$ and $P_\Ec$ denote the orthogonal projections of $H$ onto
$\clspan(\Fc)$ and $\clspan(\Ec)$, respectively, and recall that these
projections can be realized as in equation \eqref{orthogproj}.
Then for $N > N_\eps$ we have the following:
\begin{align}
& |S_N(c_N)| \, (s_N - d_N r_N) \label{reduction2} \\[1 \jot]
& \qquad \EQ \sum_{j \in S_N(c_N)} \ip{\te_j}{P_\Fc e_j} \minus
      \sum_{i \in I_N(c_N)} \ip{P_\Ec f_i}{\tf_i}
      \notag \\[1 \jot]
& \qquad \EQ \sum_{j \in S_N(c_N)} \, \sum_{i \in I} \,
      \ip{f_i}{e_j} \, \ip{\te_j}{\tf_i} \minus
      \sum_{i \in I_N(c_N)} \, \sum_{j \in J} \,
      \ip{f_i}{e_j} \, \ip{\te_j}{\tf_i} \notag \\[2 \jot]
& \qquad \EQ T_1 + T_2 - T_3 - T_4, \notag
\end{align}
where
\begin{align*}
T_1 & \EQ \sum_{j \in S_N(c_N)} \,
          \sum_{i \in I \setminus I_{N+N_\eps}(c_N)}
          \ip{f_i}{e_j} \, \ip{\te_j}{\tf_i}, \\[1 \jot]
T_2 & \EQ \sum_{j \in S_N(c_N)} \,
          \sum_{i \in I_{N+N_\eps}(c_N) \setminus I_N(c_N)}
          \ip{f_i}{e_j} \, \ip{\te_j}{\tf_i}, \allowdisplaybreaks \\[1 \jot]
T_3 & \EQ \sum_{i \in I_{N-N_\eps}(c_N)} \,
          \sum_{j \in G \setminus S_N(c_N)}
          \ip{f_i}{e_j} \, \ip{\te_j}{\tf_i}, \\[1 \jot]
T_4 & \EQ \sum_{i \in I_N(c_N) \setminus I_{N-N_\eps}(c_N)} \,
          \sum_{j \in G \setminus S_N(c_N)}
          \ip{f_i}{e_j} \, \ip{\te_j}{\tf_i}.
\end{align*}
We will estimate each of these quantities in turn.

\smallskip
\emph{Estimate $T_1$}.
If $j \in S_N(c_N)$, then $I_{N_\eps}(j) \subset I_{N+N_\eps}(c_N)$,
so by $\ell^2$-column decay we have
$$\sum_{i \in I \setminus I_{N+N_\eps}(c_N)} |\ip{f_i}{e_j}|^2
\LE \sum_{i \in I \setminus I_{N_\eps}(j)} |\ip{f_i}{e_j}|^2
\LT \eps.$$
Using this and the fact that $\set{\tf_i}_{i \in I}$ is a frame sequence,
we estimate that
\begin{align*}
|T_1|
& \LE \sum_{j \in S_N(c_N)} \,
      \biggparen{\sum_{i \in I \setminus I_{N+N_\eps}(c_N)}
                 |\ip{f_i}{e_j}|^2}^{1/2} \,
      \biggparen{\sum_{i \in I \setminus I_{N+N_\eps}(c_N)}
                 |\ip{\te_j}{\tf_i}|^2}^{1/2} \\[2 \jot]
& \LE \sum_{j \in S_N(c_N)} \, \eps^{1/2} \,
      \biggparen{\frac1A \, \norm{\te_j}^2}^{1/2}
\LE |S_N(c_N)| \, \biggparen{\frac{\eps}{AE}}^{1/2}.
\end{align*}

\smallskip
\emph{Estimate $T_2$}.
By \eqref{annulus}, we have
$|I_{N+N_\eps}(c_N) \setminus I_N(c_N)|
\le K \, \bigparen{|S_{N+N_\eps}(c_N)| - |S_N(c_N)|}$.
Since $\set{e_j}_{j \in G}$ and $\set{\te_j}_{j \in G}$ are frame
sequences, we therefore have
\begin{align*}
|T_2|
& \LE \sum_{i \in I_{N+N_\eps}(c_N) \setminus I_N(c_N)} \,
      \biggparen{\sum_{j \in G} |\ip{f_i}{e_j}|^2}^{1/2} \,
      \biggparen{\sum_{j \in G} |\ip{\te_j}{\tf_i}|^2}^{1/2} \\[2 \jot]
& \LE \sum_{i \in I_{N+N_\eps}(c_N) \setminus I_N(c_N)} \,
      \Bigparen{E \, \norm{f_i}^2}^{1/2} \,
      \Bigparen{\frac1F \, \norm{\tf_i}^2}^{1/2} \,
      \allowdisplaybreaks \\[2 \jot]
& \LE K \, \bigparen{|S_{N+N_\eps}(c_N)| - |S_N(c_N)|} \,
      \biggparen{\frac{EB}{FA}}^{1/2}.
\end{align*}

\smallskip
\emph{Estimate $T_3$}.
This estimate is similar to the one for $T_1$.
If $i \in I_{N-N_\eps}(c_N)$ then $a(i) \in S_{N-N_\eps}(c_N)$, so
$S_{N_\eps}(a(i)) \subset S_N(c_N)$.
Hence, by $\ell^2$-row decay,
$$\sum_{j \in G \setminus S_{N}(c_N)} |\ip{f_i}{e_j}|^2
\LE \sum_{j \in G \setminus S_{N_\eps}(a(i))} |\ip{f_i}{e_j}|^2
\LT \eps.$$
Since $\set{\te_j}_{j \in G}$ is a frame sequence, we therefore have
\begin{align*}
|T_3|
& \LE \sum_{i \in I_{N-N_\eps}(c_N)} \,
      \biggparen{\sum_{j \in G \setminus S_N(c_N)}
                 |\ip{f_i}{e_j}|^2}^{1/2} \,
      \biggparen{\sum_{j \in G} |\ip{\te_j}{\tf_i}|^2}^{1/2} \\[2 \jot]
& \LE \sum_{i \in I_{N-N_\eps}(c_N)} \, \eps^{1/2} \,
      \biggparen{\frac1E \, \norm{\tf_i}^2}^{1/2}
\LE  K \, |S_{N - N_\eps}(c_N)| \, \biggparen{\frac{\eps}{AE}}^{1/2}.
\end{align*}

\smallskip
\emph{Estimate $T_4$}.
This estimate is similar to the one for $T_2$.
Since $\set{e_j}_{j \in G}$ and $\set{\te_j}_{j \in G}$ are frame
sequences, we have for $N > N_\eps$ that
\begin{align*}
|T_4|
& \LE \sum_{i \in I_N(c_N) \setminus I_{N-N_\eps}(c_N)} \,
      \biggparen{\sum_{j \in G} |\ip{f_i}{e_j}|^2}^{1/2} \,
      \biggparen{\sum_{j \in G} |\ip{\te_j}{\tf_i}|^2}^{1/2} \\[2 \jot]
& \LE K \, \bigparen{|S_N(c_N)| - |S_{N-N_\eps}(c_N)|} \,
      \biggparen{\frac{EB}{FA}}^{1/2}.
\end{align*}

\smallskip
\emph{Final Estimate}.
Applying the above estimates to \eqref{reduction2},
we find that if $N > N_\eps$, then
\begin{align*}
|s_N - d_N r_N|
& \LE \frac{|T_1| + |T_2| + |T_3| + |T_4|} {|S_N(c_N)|} 
\allowdisplaybreaks \\[2 \jot]
& \LE \biggparen{\frac{\eps}{AE}}^{1/2} \plus
      K \, \biggparen{\frac{EB}{FA}}^{1/2} \,
         \frac{|S_{N+N_\eps}(c_N)| - |S_N(c_N)|}{|S_N(c_N)|} \, \plus
      \\[1 \jot]
& \qquad\qquad
      K \, \biggparen{\frac{\eps}{AE}}^{1/2} \,
         \frac{|S_{N-N_\eps}(c_N)|}{|S_N(c_N)|} \, \plus
      \\[1 \jot]
& \qquad\qquad
      K \, \biggparen{\frac{EB}{FA}}^{1/2} \,
         \frac{|S_N(c_N)| - |S_{N-N_\eps}(c_N)|}{|S_N(c_N)|}.
\end{align*}
Consequently, applying the asymptotics in \eqref{asymptotics},
we conclude that
$$\limsup_{N \to \infty} |s_N - d_N r_N|
\LE \biggparen{\frac{\eps}{AE}}^{1/2} \plus 0 \plus
    K \biggparen{\frac{\eps}{AE}}^{1/2} \plus 0.$$
Since $\eps$ was arbitrary, this implies
$\lim_{N \to \infty} (s_N - d_N r_N) = 0$, as desired.

\medskip
(b) Since ultrafilter limits exist for any bounded sequence and furthermore
are linear and respect products, we have
\begin{align*}
0 \EQ \plim_{N \in \N} (s_N - d_N r_N)
& \EQ \Bigparen{\plim_{N \in \N} s_N} \minus
      \Bigparen{\plim_{N \in \N} d_N} \,
      \Bigparen{\plim_{N \in \N} r_N} \\[1 \jot]
& \EQ \cM_\Fc(\Ec;p,c) \minus D(p,c) \cdot \cM_\Ec(\Fc;p,c).
      \qedhere
\end{align*}
\end{proof}

\subsection{Applications of the Density--Relative Measure Theorem}
\label{applications}

In this section we will derive some consequences of
Theorem~\ref{densityredundancy}.

Our first result specializes Theorem~\ref{densityredundancy} to the case
where $\Fc$ and $\Ec$ are both frames for~$H$, including the important
special cases where $\Ec$ is actually a Riesz basis for~$H$.
It also connects the infinite excess result of
Proposition~\ref{infiniteexcess}.

\begin{theorem}[Abstract Density Theorem] \label{excesscor}
Let $\Fc = \set{f_i}_{i \in I}$ and $\Ec = \set{e_j}_{j \in G}$ be frames
for~$H$, and let $a \colon I \to G$ be an associated map such that
$D^+(I,a) < \infty$.
If $(\Fc, a, \Ec)$ has both $\ell^2$-column decay and $\ell^2$-row decay,
then the following statements hold.

\smallskip
\begin{enumerate}
\item[(a)]
For each free ultrafilter $p$ and
sequence of centers $c = (c_N)_{N \in \N}$ in $G$, we have
\begin{equation} \label{equality1}
\cM(\Ec;p,c) \EQ D(p,c)   \cdot \cM(\Fc;p,c).
\end{equation}
Consequently,
\begin{align}
\cM^-(\Ec) & \LE D^+(I,a) \cdot \cM^-(\Fc) \LE \cM^+(\Ec), \label{equality2}
             \\[2 \jot]
\cM^-(\Ec) & \LE D^-(I,a) \cdot \cM^+(\Fc) \LE \cM^+(\Ec). \label{equality3}
\end{align}

\item[(b)]
If $D^+(I,a) > \cM^+(\Ec)$, then there exists an infinite set $J \subset I$
such that $\set{f_i}_{i \in I \setminus J}$ is still a frame for $H$.
\end{enumerate}

\medskip\noindent
If $\Ec$ is a Riesz basis for $H$ then the following additional
statements hold.

\smallskip
\begin{enumerate}
\item[(c)]
For each free ultrafilter $p$ and
sequence of centers $c = (c_N)_{N \in \N}$ in $G$, we have
$$\cM(\Fc;p,c) \EQ \frac1{D(p,c)}, \quad
\cM^-(\Fc) \EQ \frac1{D^+(I,a)}, \quad
\cM^+(\Fc) \EQ \frac1{D^-(I,a)}.$$

\medskip
\item[(d)] $D^-(I,a) \ge 1$.

\medskip
\item[(e)] 
If $D^+(I,a) > 1$, then there exists an infinite subset
$J \subset I$ such that $\set{f_i}_{i \in I \setminus J}$
is still a frame for~$H$.

\medskip
\item[(f)] 
If $\Fc$ is also a Riesz basis for $H$,
then for each free ultrafilter $p$ and
sequence of centers $c = (c_N)_{N \in \N}$ in $G$, we have
\begin{align*}
& D^-(I,a) \EQ D(p,c) \EQ D^+(I,a) \EQ 1, \\[1 \jot]
& \cM^-(\Fc) \EQ \cM(\Fc;p,c) \EQ \cM^+(\Fc) \EQ 1.
\end{align*}
\end{enumerate}
\end{theorem}
\begin{proof}
(a) Since the closed span of $\Fc$ and $\Ec$ is all of $H$, the
equality in \eqref{equality1}
is a restatement of Theorem~\ref{densityredundancy}(a).
For the first inequality in \eqref{equality2}, choose an ultrafilter
$p$ and sequence of centers $c$ such that
$\cM^-(\Fc) = \cM(\Fc;p,c)$.
Then we have
$$\cM^-(\Ec)
\LE \cM(\Ec;p,c)
\EQ D(p,c) \cdot \cM(\Fc;p,c)
\LE D^+(p,c) \cdot \cM^-(\Fc).$$
The other inequalities in \eqref{equality2} and \eqref{equality3}
are similar.

\medskip
(b) In this case it follows from \eqref{equality2} that
$\cM^-(\Fc) \le \cM^+(\Ec) / D^+(I,a) < 1$,
so the result follows from Proposition~\ref{infiniteexcess}.

\medskip
(c) If $\Ec$ is a Riesz basis then
$\cM(\Ec;p,c) = \cM^\pm(\Ec) = 1$,
so the result follows from part~(a).

\medskip
(d) Follows from part~(c) and the fact that
$0 \le \cM^+(\Fc) \le 1$.

\medskip
(e) Follows from part~(b) and the fact that
$\cM^+(\Ec) = 1$.

\medskip
(f) If $\Fc$ is a Riesz basis then $\cM^\pm(\Fc) = 1$, so this
follows from part~(c).
\end{proof}

Note that the conclusion $D^-(I,a) \ge 1$ of Theorem~\ref{excesscor}(d)
is shown under a weaker hypothesis in Theorem~\ref{necessary}.
Specifically, Theorem~\ref{necessary} requires only the hypothesis that the
weak HAP be satisfied.
However, the stronger localization hypotheses of Theorem~\ref{excesscor}
($\ell^2$-column and row decay) yields the significantly
stronger conclusions of Theorem~\ref{excesscor}.

Next we derive relationships among the density, frame bounds, and
norms of the frame elements for localized frames.
In particular, part~(a) provides an estimate of the relations between frame
bounds, density, and limits of averages of the norms of frame elements.
Many of the frames that are important in applications, such as
Gabor frames, are uniform norm frames, i.e., all the frame elements
have identical norms, and for these frames these averages are a constant.
As a consequence, we show that if $\Fc$ and $\Ec$ are both tight uniform
norm frames, then the index set $I$ must have uniform density.

\begin{theorem}[Density--Frame Bounds] \label{framebounds}
Let $\Fc = \set{f_i}_{i \in I}$ be a frame for $H$ with frame bounds $A$,~$B$,
and let $\Ec = \set{e_j}_{j \in G}$ be a frame for $H$ with frame bounds
$E$,~$F$.
Let $a \colon I \to G$ be an associated map such that $D^+(I,a) < \infty$.
If $(\Fc, a, \Ec)$ has both $\ell^2$-column decay and $\ell^2$-row decay,
then the following statements hold.

\smallskip
\begin{enumerate}
\item[(a)]
For each free ultrafilter $p$ and
sequence of centers $c = (c_N)_{N \in \N}$ in $G$, we have
\begin{align}
\frac1{F} \,
\plim_{N \in \N} \frac1{|S_N(c_N)|} \sum_{j \in S_N(c_N)} \norm{e_j}^2
& \LE \frac{D(p,c)}{A} \,
      \plim_{N \in \N} \frac1{|I_N(c_N)|} \sum_{i \in I_N(c_N)} \norm{f_i}^2,
      \label{firstplim} \\[1 \jot]
\frac1{E} \,
\plim_{N \in \N} \frac1{|S_N(c_N)|} \sum_{j \in S_N(c_N)} \norm{e_j}^2
& \GE \frac{D(p,c)}{B} \,
      \plim_{N \in \N} \frac1{|I_N(c_N)|} \sum_{i \in I_N(c_N)} \norm{f_i}^2.
      \label{secondplim}
\end{align}

\smallskip
\item[(b)]
We have
$$\frac{A}{F} \, \frac{\liminf_j \norm{e_j}^2}{\limsup_i \norm{f_i}^2}
\LE D^-(I,a)
\LE D^+(I,a)
\LE \frac{B}{E} \, \frac{\limsup_j \norm{e_j}^2}{\liminf_i \norm{f_i}^2}.$$

\smallskip
\item[(c)]
If $\Fc$ and $\Ec$ are both uniform norm frames, with
$\norm{f_i}^2 = \mathcal{N}_\Fc$ for $i \in I$ and
$\norm{e_j}^2 = \mathcal{N}_\Ec$ for $j \in G$, then
$$\frac{A \, \mathcal{N}_\Ec}{F \, \mathcal{N}_\Fc}
\LE D^-(I,a)
\LE D^+(I,a)
\LE \frac{B \, \mathcal{N}_\Ec}{E \, \mathcal{N}_\Fc}.$$
Consequently, if $\Fc$ and $\Ec$ are both tight uniform norm frames,
then $I$ has uniform density, with
$D^-(I,a) = D^+(I,a) = (A \, \mathcal{N}_\Ec)/(E \, \mathcal{N}_\Fc)$.
\end{enumerate}
\end{theorem}

\begin{proof}
(a) Let $S$ be the frame operator for $\Fc$.
Then $A\one \le S \le B\one$, so we have
$\ip{f_i}{\tf_i}
= \ip{f_i}{S^{-1}(f_i)}
\le \frac1{A} \, \ip{f_i}{f_i}
= \frac1{A} \, \norm{f_i}^2$,
and hence
$$\cM(\Fc;p,c)
\EQ \plim_{N \in \N} \frac1{|I_N(c_N)|} \sum_{i \in I_N(c_N)} \ip{f_i}{\tf_i}
\LE \frac1{A} \,
    \plim_{N \in \N} \frac1{|I_N(c_N)|} \sum_{i \in I_N(c_N)} \norm{f_i}^2.$$
Similarly $\ip{\te_j}{e_j} \ge \frac1{F} \, \norm{e_j}^2$, so
$$\cM(\Ec;p,c)
\EQ \plim_{N \in \N} \frac1{|S_N(c_N)|} \sum_{j \in S_N(c_N)} \ip{\te_j}{e_j}
\GE \frac1{F} \,
    \plim_{N \in \N} \frac1{|S_N(c_N)|} \sum_{j \in S_N(c_N)} \norm{e_j}^2.$$
Combining these inequalities with the equality
$\cM(\Ec;p,c) = D(p,c) \cdot \cM(\Fc;p,c)$ from Theorem~\ref{excesscor}(a)
yields \eqref{firstplim}.
Inequality \eqref{secondplim} is similar, using
$\ip{f_i}{\tf_i} \ge \frac1B \, \norm{f_i}^2$ and
$\ip{\te_j}{e_j} \le \frac1E \, \norm{e_j}^2$.

\medskip
(b) Observe that
$$\liminf_{i \in I} \, \norm{f_i}^2
\LE \plim_{N \in \N} \frac1{|I_N(c_N)|} \sum_{i \in I_N(c_N)} \norm{f_i}^2
\LE \limsup_{i \in I} \, \norm{f_i}^2,$$
and combine this and a similar inequality for $\Ec$ with \eqref{firstplim}.

\medskip
(c) This is an immediate consequence of part~(b).
\end{proof}

A similar result can be formulated in terms of the norms $\norm{\tf_i}$
of the canonical dual frame elements, by using the inequality
$A \, \norm{\tf_i}^2 \le \ip{f_i}{\tf_i} \le B \, \norm{\tf_i}$.

\subsection{Removing Sets of Positive Measure} \label{positivesec}

In this section, we will show that by imposing a stronger form of
localization than we used in Theorem~\ref{excesscor},
a subset of positive measure may be removed yet still leave a frame.
This is a stronger conclusion than the infinite excess statements
of Proposition~\ref{infiniteexcess} or Theorem~\ref{excesscor},
which only state that an infinite set may be removed, without any
conclusion about the density of that set.

In the remainder of this section we will use the results of
Appendix~\ref{selflocappend}, as well as the following notations.
If $\Fc = \set{f_i}_{i \in I}$ is a frame then
the orthogonal projection of $\ell^2(I)$ onto the range of the analysis
operator $T$ is $\Pb = T S^{-1} T^*$.
Given $J \subset I$, we define truncated analysis and frame operators
$T_J f = \set{\ip{f}{f_i}}_{i \in J}$ and
$S_J f = \sum_{i \in J} \ip{f}{f_i} \, f_i$.
We let $R_J \colon \ell^2(I) \to \ell^2(I)$ be the projection operator
given by $(R_J c)_k = c_k$ for $k \in J$, and $0$ otherwise.
Written as matrices,
$$\Pb \EQ T S^{-1} T^* \EQ [\ip{f_i}{\tf_j}]_{i,j \in I}
\qquad\text{and}\qquad
T_J S^{-1} T_J^* \EQ [\ip{f_i}{\tf_j}]_{i,j \in J}.$$

The following lemma characterizes those subsets of a frame which
can be removed yet still leave a frame.

\begin{lemma} \label{frameremove}
Let $\Fc = \set{f_i}_{i \in I}$ be a frame for $H$,
with frame bounds $A$, $B$.
Let $J \subset I$ be given, and define
\begin{equation} \label{rhodef}
\rho
\EQ \norm{T_J S^{-1} T_J^*}
\EQ \norm{S^{-1/2} S_J S^{-1/2}}
\EQ \norm{R_J \Pb R_J}.
\end{equation}
Then $\Fc_{I \setminus J} = \set{f_i}_{i \in I \setminus J}$
is a frame for $H$ if and only if $\rho < 1$.
In this case, $A(1-\rho)$,~$B$ are frame bounds for $\Fc_{I \setminus J}$.
\end{lemma}
\begin{proof}
First, the fact that equality holds in \eqref{rhodef} is a consequence
of the fact that $\norm{L^* L} = \norm{LL^*}$ for any operator~$L$.
Specifically,
\begin{align*}
\norm{S^{-1/2} S_J S^{-1/2}}
& \EQ \norm{(S^{-1/2} T_J^*) (S^{-1/2} T_J^*)^*}
  \EQ \norm{(S^{-1/2} T_J^*)^* (S^{-1/2} T_J^*)}
      \allowdisplaybreaks \\[1 \jot]
& \EQ \norm{T_J S^{-1} T_J^*}
  \EQ \norm{R_J T S^{-1} T^* R_J}
  \EQ \norm{R_J \Pb R_J}.
\end{align*}
Second,
since $\Fc_{I \setminus J}$ is a subset of $\Fc$, it is clearly a
Bessel sequence with Bessel bound $B$.
Further, $S_{I \setminus J}$ is a bounded operator on $H$, satisfying
$0 \le S_{I \setminus J} \le S \le BI$.
Therefore, $\Fc_{I \setminus J}$ is a frame for $H$ with frame bounds
$A'$, $B$ if and only if $A' \one \le S_{I \setminus J}$.

Suppose now that $\rho = \norm{S^{-1/2} S_J S^{-1/2}} < 1$.
Then
$$S_{I \setminus J}
\EQ S - S_J
\EQ S^{1/2} (\one - S^{-1/2} S_J S^{-1/2}) S^{1/2}$$
is invertible.
Further,
$$\ip{S^{-1/2} S_J S^{-1/2}f}{f}
\LE \norm{S^{-1/2} S_J S^{-1/2}} \, \norm{f}^2
\LE \rho \, \norm{f}^2
\EQ \ip{\rho\one f}{f},$$
so
\begin{align*}
S_{I \setminus J}
& \EQ S^{1/2} (\one - S^{-1/2} S_J S^{-1/2}) S^{1/2} \\
& \GE S^{1/2} (\one - \rho\one) S^{1/2}
\EQ (1-\rho) S
\GE (1-\rho) A\one.
\end{align*}
Thus $\Fc_{I \setminus J}$ is a frame for $H$ with frame bounds
$(1-\rho)A$, $B$.

Conversely, if $\Fc_{I \setminus J}$ is a frame with frame bounds
$A'$, $B$ then $S_{I \setminus J} \ge A'\one$, so
$$\one - S^{-1/2} S_J S^{-1/2}
\EQ S^{-1/2} S_{I \setminus J} S^{-1/2}
\GE S^{-1/2} A' \one S^{-1/2}
\EQ A' S^{-1}
\GE \frac{A'}B \one.$$
Hence
$\rho
= \norm{S^{-1/2} S_J S^{-1/2}}
\le \norm{(1 - \frac{A'}B) \one}
= 1 - \frac{A'}B
< 1$.
\end{proof}

Now we can give the first main result of this section, that if
$\cM(\Fc^+) < 1$ and we have $\ell^1$-localization with respect to the
dual frame, then a set of positive uniform density can be removed
yet still leave a frame.
Note by Theorem~\ref{selflocthm} the hypothesis of $\ell^1$-localization with 
respect to the canonical dual is implied by $\ell^1$-self-localization.
Although we omit it, it is possible to give a direct proof
of the following result under the hypothesis of $\ell^1$-self-localization
that does not appeal to Theorem~\ref{selflocthm}.

\begin{theorem}[Positive Uniform Density Removal] \label{positiveremove}
Let $\Fc = \set{f_i}_{i \in I}$ be a frame sequence with frame bounds $A$,~$B$,
with associated map $a \colon I \to G$,
and assume that the following statements hold:

\smallskip
\begin{enumerate}
\item[(a)]
$0 < D^-(I,a) \le D^+(I,a) < \infty$,

\medskip
\item[(b)]
$(\Fc,a)$ is $\ell^1$-localized with respect to its canonical dual frame, and

\medskip
\item[(c)]
$\cM^+(\Fc) < 1$.
\end{enumerate}

\smallskip\noindent
Then there exists a subset $J \subset I$ such that
$D^+(J,a) = D^-(J,a) > 0$ and
$\Fc_{I \setminus J} = \set{f_i}_{i \in I \setminus J}$
is a frame for $\clspan(\Fc)$.

Moreover, if $\cM^+(\Fc) < \alpha < 1$ and $J_\alpha$ is defined
by~\eqref{Jalphadef}, i.e.,
$$J_\alpha \EQ \set{i \in I : \ip{f_i}{\tf_i} \le \alpha},$$ 
then for each $0 < \eps < 1-\alpha$
there exists a subset $J \subset J_\alpha$ such that
$D^+(J,a) = D^-(J,a) > 0$ and
$\Fc_{I \setminus J} = \set{f_i}_{i \in I \setminus J}$
is a frame for $\clspan(\Fc)$ with frame bounds $A(1-\alpha-\eps)$,~$B$.
\end{theorem}

\begin{proof}

Note first that by Corollary~\ref{Jalphacoro}(a),
if $\cM^+(\Fc) < \alpha < 1$ then we have that $D^-(J_\alpha,a) > 0$.
Also, since $(\Fc,a)$ is $\ell^1$-localized with respect to its dual frame,
there exists $r \in \ell^1(G)$ such that
$|\ip{f_i}{\tf_j}| \le r_{a(i) - a(j)}$ for all $i$, $j \in I$.
Given $0 < \eps < 1 - \alpha$, let $N_\eps$ be large enough that
$$\sum_{k \in G \setminus S_{N_\eps}(0)} r_k
\LT \eps.$$
Since $D^-(J_\alpha,a) > 0$, there exists $N_0 > 0$ such that
$|I_{N_0}(j) \cap J_\alpha| > 0$ for every $j \in G$.
Let $N = \max\set{N_\eps,N_0}$, and define
$$\Qc = \set{S_N(2Nk) : k \in G}.$$
Each preimage $I_N(2Nk) = a^{-1}(S_N(2Nk))$ of the boxes in $\Qc$
contains at least one point of~$J_\alpha$.
For each $k$, select one such point, say~$i_k \in I_N(2Nk) \cap J_\alpha$,
and set
$J = \set{i_k : k \in G}$.
Then $J$ has positive density, with
$D^+(J,a) = D^-(J,a) = \frac1{|S_{2N}(0)|}$.

Consider now the matrix
$T_J S^{-1} T_J^* = [\ip{f_i}{\tf_j}]_{i,j \in J}$.
Write $T_J S^{-1} T_J^* = D + V$, where $D$ is the diagonal part
of $T_J S^{-1} T_J^*$ and $V=[v_{ij}]_{i,j\in J}$.
By the definition of $J_\alpha$, we have
$\norm{D} = \sup_{i \in J} \ip{f_i}{\tf_i} \le \alpha$.
Define
$$s_k \EQ
\begin{cases}
r_k, & k \notin S_{N_\eps}(0), \\
0,   & k \in S_{N_\eps}(0).
\end{cases}$$
If $i$, $j \in J$ and $i \ne j$, then $a(i) - a(j) \notin S_{N_\eps}(0)$,
and therefore
$|v_{ij}| = |\ip{f_i}{\tf_j}| \le r_{a(i) - a(j)} = s_{a(i) - a(j)}$.
On the other hand, $|v_{ii}| = 0 = s_{a(i) - a(i)}$.
Applying Proposition~\ref{matrixdecay}(a) to $V$ and the index set $J$
therefore yields
$$\norm{V}
\LE \sum_{k \in G} s_k
\EQ \sum_{k \in G \setminus S_{N_\eps}(0)} r_k
\LT \eps.$$
Therefore
$\norm{T_J S^{-1} T_J^*}
\le \norm{D} + \norm{V}
\le \alpha + \eps
< 1$.
Lemma~\ref{frameremove} therefore implies that
$\set{f_i}_{i \in I \setminus J}$ is a frame for $H$
with frame bounds $A(1-\alpha-\eps)$, $B$.
\end{proof}

If we impose $\ell^2$-column decay and $\ell^2$-row decay,
then we can reformulate Theorem~\ref{positiveremove}
in terms of density instead of relative measure.

\begin{corollary} \label{removecorollary}
Let $\Fc = \set{f_i}_{i \in I}$ and $\Ec = \set{e_j}_{j \in G}$ be frames
for $H$, and let $A$,~$B$ be frame bounds for $\Fc$.
Let $a \colon I \to G$ be an associated map, and
assume that the following statements hold:

\smallskip
\begin{enumerate}
\item[(a)]
$0 < D^-(I,a) \le D^+(I,a) < \infty$,

\medskip
\item[(b)]
$(\Fc,a)$ is $\ell^1$-localized with respect to its canonical dual frame,

\medskip
\item[(c)]
$(\Fc,a,\Ec)$ has both $\ell^2$-column decay and $\ell^2$-row decay, and

\medskip
\item[(d)]
$\cM^+(\Ec) < D^-(I,a)$;
in particular, $D^-(I,a) > 1$ if $\Ec$ is a Riesz basis.
\end{enumerate}

\smallskip\noindent
Then $\cM^+(\Fc) < 1$, and
then there exists a subset $J \subset I$ such that
$D^+(J,a) = D^-(J,a) > 0$ and
$\Fc_{I \setminus J} = \set{f_i}_{i \in I \setminus J}$
is a frame for $\clspan(\Fc)$.

Moreover, if $\cM^+(\Fc) < \alpha < 1$ and $J_\alpha$ is defined
by~\eqref{Jalphadef}, then for each $0 < \eps < 1-\alpha$
there exists a subset $J \subset J_\alpha$ such that
$D^+(J,a) = D^-(J,a) > 0$ and
$\Fc_{I \setminus J} = \set{f_i}_{i \in I \setminus J}$
is a frame for $\clspan(\Fc)$ with frame bounds $A(1-\alpha-\eps)$,~$B$.
\end{corollary}
\begin{proof}
By Theorem~\ref{excesscor} we have 
$\cM^+(\Fc) \le \frac{\cM^+(\Ec)}{D^-(I,a)} < 1$,
so the result follows by applying Theorem~\ref{positiveremove}.
\end{proof}

Theorem~\ref{positiveremove} and Corollary~\ref{removecorollary} are
evidence that the reciprocal of the relative measure should in fact be
a quantification of the redundancy of an abstract frame.
Concentrating for purposes of discussion on the case where $\Ec$ is a
Riesz basis (and hence $\cM^+(\Ec) = 1$),
this quantification would be precise if it was the case that if
$\Fc = \set{f_i}_{i \in I}$ is an appropriately localized frame and if
$\cM^+(\Fc) < 1$, then there exists a subset $I'$ of $I$ with density
$1+\eps$ such that $\Fc' = \set{f_i}_{i \in I'}$ is still a frame for $H$
(and not merely, as implied by Theorem~\ref{positiveremove} or
Corollary~\ref{removecorollary}, that there is some set $J$ with positive
density such that $\set{f_i}_{i \in I \setminus J}$ is a frame).
To try to prove such a result, we could attempt to iteratively apply
Corollary~\ref{removecorollary}, repeatedly removing sets of positive
measure until we are left with a subset of density $1+\eps$ that is
still a frame.
However there are several obstructions to this approach.
One is that with each iteration, the lower frame bound is reduced and
may approach zero in the limit.
A second problem is that the lower density of $I'$ may eventually approach~$1$.
Because Corollary~\ref{removecorollary} removes sets of uniform density,
we would then have $D^+(I',a)$ approaching $1 + D^+(I,a)- D^-(I,a)$,
which for a frame with non-uniform density would not be of the form
$1+\eps$ with $\eps$ small.
Due to the length and breadth of this work, we have chosen to omit some
results dealing with this second obstruction.

\subsection{Localized Frames and $\eps$-Riesz sequences} \label{rieszsec}

Feichtinger has conjectured that every frame that is norm-bounded below
can be written as a union of a finite number of Riesz sequences
(systems that are Riesz bases for their closed linear spans).
It is shown in \cite{CCLV03}, \cite{CV03}, \cite{CT05}
that Feichtinger's conjecture equivalent to the
celebrated Kadison--Singer (paving) conjecture.
and that both of these are equivalent to a conjectured generalization 
of the Bourgain--Tzafriri restricted invertibility theorem.

In this section we will show that every $\ell^1$-self-localized frame
that is norm-bounded below is a finite union of $\eps$-Riesz sequences,
and every frame that is norm-bounded below and $\ell^1$-localized with
respect to its dual frame is a finite union of Riesz sequences.

\begin{definition}
If $0 < \eps < 1$ and $f_i \in H$, then $\set{f_i}_{i \in I}$
is an \emph{$\eps$-Riesz sequence}
if there exists a constant $A > 0$ such that for every sequence
$(c_i)_{i \in I} \in \ell^2(I)$ we have
$$(1-\eps)A \, \sum_{i \in I} |c_i|^2
\LE \Bignorm{\sum_{i \in I} c_i f_i}^2
\LE (1+\eps)A \, \sum_{i \in I} |c_i|^2.
\qquad\qed$$
\end{definition}

\smallskip
Every $\eps$-Riesz sequence is a Riesz sequence, i.e.,
a Riesz basis for its closed linear span.

\begin{theorem} \label{epsrieszthm}
Let $\Fc = \set{f_i}_{i \in I}$ be a sequence in $H$ and let
$a \colon I \to G$ be an associated map.
If

\begin{enumerate}
\item[(a)] $(\Fc,a)$ is $\ell^1$-self-localized,

\smallskip
\item[(b)] $D^+(I,a) < \infty$, and

\smallskip
\item[(c)] $\inf_i \norm{f_i} > 0$,
\end{enumerate}

\smallskip\noindent
then for each $0 < \eps < \inf_i \norm{f_i}$,
$\Fc$ can be written as a finite union of $\eps$-Riesz sequences.
\end{theorem}
\begin{proof}
Recall that $G$ has the form
$G \EQ \prod_{i=1}^d a_i \Z \, \times \, \prod_{j=1}^e \Z_{n_j}$.
For simplicity of notation, we will treat the case where $a_i = 1$
for all $i$, so $G = \Z^d \times H$ with 
$H = \prod_{j=1}^e \Z_{n_j}$.
The general case is similar.

For this proof we will use boxes in $G$ of the form
$$B_N(j)
\EQ j \plus \biggparen{\biggl[-\frac{N}{2},\frac{N}{2}\biggr)^d
                       \times H},
\qquad j \in G, \ N > 0.$$

Set $m = \inf_i \norm{f_i}^2$ and
$M = \sup_i \norm{f_i}^2$.
Fix $0 < \eps < m$, set $\delta = \eps m$, and choose $K$ so that
$\frac{M-m}K < \frac{\delta}2$.
Partition $I$ into subsequences $\set{J_k}_{k=1}^K$ so that
$$\forall\, i \in J_k, \qquad
m + \frac{M-m}K \, (k-1)
\LE \norm{f_i}^2
\LE m + \frac{M-m}K \, k.$$

Since $(\Fc,a)$ is $\ell^1$-self-localized, there exists an
$r \in \ell^1(G)$ such that
$|\ip{f_i}{f_j}| \le r_{a(i) - a(j)}$
for all $i$, $j \in I$.
Let $N_\delta$ be large enough that
$$\sum_{n \in G \setminus B_{N_\delta}(0)} r_n
\LT \frac{\delta}2.$$

Let $\set{u_\nu}_{\nu=1}^{2^d}$ be a list of the vertices of
the unit cube $[0,1]^d$.
For $\nu = 1, \dots, 2^d$, define
$$\Qc_\nu
\EQ \set{B_{N_\delta}(2N_\delta n + N_\delta u_\nu)}_{n \in \Z^d}.$$
Each $\Qc_\nu$ is a set of disjoint boxes in $G$, each of which is separated
by a distance of at least~$N_\delta$ from the other boxes.
Furthermore, the union of the boxes in $\Qc_\nu$ for
$\nu = 1,\ldots,2^d$ forms a disjoint cover of $G$. 

Since $D^+(I,a) < \infty$, we have
$L = \sup_{n \in G} |I_{N_\delta}(n)| < \infty$. 
Therefore each box in $\Qc_\nu$ contains at most $L$ points of $a(I)$.
By choosing, for each fixed $k$ and $\nu$, at most a single element of 
$J_k$ out of each box in $\Qc_\nu$, we can divide each subsequence 
$J_k$ into $2^d L$ or fewer subsequences
$\set{J_{k\ell}}_{\ell=1}^{K_k}$ in such a way that
$$\forall\, i, j \in J_{k\ell}, \qquad
i \ne j \ \implies\ a(i) - a(j) \notin B_{N_\delta}(0).$$

Fix $k$, $\ell$,
let $G_{k\ell} = [\ip{f_i}{f_j}]_{i,j \in J_{k\ell}}$, and write
$G_{k\ell} = D_{k\ell} + V_{k\ell}$, where $D_{k\ell}$ is the diagonal
part of~$G_{k\ell}$.
Set
$$s_n \EQ
\begin{cases}
r_n, & n \notin B_{N_\delta}(0), \\
0,   & n \in B_{N_\delta}(0),
\end{cases}$$
If we write the entries of $V_{k\ell}$ as $V_{k\ell} = [v_{ij}]_{i,j \in J}$
then we have $|v_{ij}| \le s_{a(i) - a(j)}$ for all $i$, $j \in J$.
Applying Proposition~\ref{matrixdecay} to the matrix $V_{k\ell}$ and the
index set $J$ therefore implies
$$\norm{V_{k\ell}}
\LE \sum_{n \in G} s_n
\EQ \sum_{n \in G \setminus B_{N_\delta}(0)} r_n 
\LT \frac{\delta}2.$$
Hence, given any sequence
$c = (c_i)_{i \in J_{k\ell}} \in \ell^2(J_{k\ell})$, we have
\begin{align*}
\biggnorm{\sum_{i \in J_{k\ell}} c_i f_i}^2
& \EQ \biggip{\sum_{i \in J_{k\ell}} c_i f_i} {\sum_{j \in J_{k\ell}} c_j f_j}
      \\[1 \jot]
& \EQ \sum_{i \in J_{k\ell}} |c_i|^2 \, \norm{f_i}^2 \plus
      \sum_{i, j \in J_{k\ell}, \, i \ne j} c_i \bar{c}_j \, \ip{f_i}{f_j}
      \allowdisplaybreaks \\[1 \jot]
& \LE \Bigparen{m + \frac{M-m}K \, k} \sum_{i \in J_{k\ell}} |c_i|^2
      \plus \ip{V_{k\ell}c}c
      \allowdisplaybreaks \\[1 \jot]
& \LE \Bigparen{m + \frac{M-m}K \, k + \frac{\delta}2} \, \norm{c}_{\ell^2}^2
      \allowdisplaybreaks \\[1 \jot]
& \LE \Bigparen{m + \frac{M-m}K \, k + \eps m} \, \norm{c}_{\ell^2}^2
      \\[1 \jot]
& \LE (1 + \eps) \, \Bigparen{m + \frac{M-m}K} \, \norm{c}_{\ell^2}^2.
\end{align*}
Similarly,
\begin{align*}
\biggnorm{\sum_{i \in J_{k\ell}} c_i f_i}^2
& \GE \Bigparen{m + \frac{M-m}K \, (k-1)} \sum_{i \in J_{k\ell}} |c_i|^2
      \minus \ip{V_{k\ell}c}c \\
& \GE \Bigparen{m + \frac{M-m}K \, k - \frac{M-m}K - \frac{\delta}2} \,
      \norm{c}_{\ell^2}^2
      \allowdisplaybreaks \\[1 \jot]
& \GE \Bigparen{m + \frac{M-m}K \, k - \delta} \, \norm{c}_{\ell^2}^2
      \allowdisplaybreaks \\[1 \jot]
& \GE \Bigparen{m + \frac{M-m}K \, k - \eps m} \, \norm{c}_{\ell^2}^2
      \\[1 \jot]
& \GE (1 - \eps) \, \Bigparen{m + \frac{M-m}K \, k} \, \norm{c}_{\ell^2}^2.
\end{align*}
Thus each $\set{f_i}_{i \in J_{k\ell}}$ is an $\eps$-Riesz sequence.
\end{proof}

\begin{corollary}
Let $\Fc = \set{f_i}_{i \in I}$ be a sequence in $H$ and let
$a \colon I \to G$ be an associated map.
If

\begin{enumerate}
\item[(a)]
$(\Fc,a)$ is $\ell^1$-localized with respect to its canonical dual frame,

\smallskip
\item[(b)] $D^+(I,a) < \infty$, and

\smallskip
\item[(c)] $\inf_i \norm{f_i} > 0$,
\end{enumerate}

\smallskip\noindent
then $\Fc$ can be written as a finite union of Riesz sequences.
\end{corollary}
\begin{proof}
Let $S$ be the frame operator for $\Fc$.
Then $(S^{-1/2}(\Fc),a)$ is $\ell^1$-self-localized by
Remark~\ref{selfremark}(b), and we have
$\inf_i \norm{S^{-1/2} f_i} > 0$ since $S^{-1/2}$ is a
continuous bijection.
If we fix $0 < \eps < \inf_i \norm{S^{-1/2}(f_i)}^2$, then
Theorem~\ref{epsrieszthm} implies that $S^{-1/2}(\Fc)$ is a finite union of
$\eps$-Riesz sequences, and hence~$\Fc$ is a finite union of Riesz sequences.
\end{proof}

\appendix
\section{The Algebra of $\ell^1$-Localized Operators}
\label{selflocappend}

Our goal in this appendix is to prove Theorem~\ref{selflocthm}.
However, we first develop some machinery about the algebra of matrices
which are bounded by Toeplitz-like matrices which have an $\ell^1$-decay
on the diagonal.

\begin{definition}
Let $I$ be a countable index set and $a \colon I \to G$ an
associated map.
We say that an $I \times I$ matrix
$V = [v_{ij}]_{i,j \in J}$ has \emph{$\ell^1$-decay}
if there exists $r \in \ell^1(G)$ such that
$|v_{ij}| \le r_{a(i) - a(j)}$.
We call $r$ an \emph{associated sequence}.
We define
$$\Bc_1(I,a) \EQ \set{V : V \text{ has } \ell^1\text{-decay}}.$$
and set $\Bc_1(G) = \Bc_1(G,Id)$, where $Id$ is the identity map.
\qed
\end{definition}

\begin{remark}
Let $\Fc = \set{f_i}_{i \in I}$ be a frame for $H$.
Let $T$ be the analysis operator and $S = T^* T$ the frame operator,
and $\tFc = \set{\tf_i}_{i \in I}$ the canonical dual frame.

\medskip
(a) $(\Fc,a)$ is $\ell^1$-self-localized
if and only if its Gram operator $V = T T^* = [\ip{f_i}{f_j}]_{i,j \in I}$
lies in~$\Bc_1(I,a)$.

\medskip
(b) The Gram operator of $\tFc$ is
$\tilde{V} = [\ip{\tf_i}{\tf_j}]_{i,j \in I} = T S^{-2} T^*$.
Since $V \tilde{V} = T S^{-1} T^* = P_V$,
the orthogonal projection onto the range of $V$,
we have that $\tilde{V} = V^\dagger$ is the pseudo-inverse of $V$.

\medskip
(c) $(\Fc,a)$ is $\ell^1$-localized with respect to its
canonical dual frame $\tFc$ if and only if the cross-Grammian matrix
$P_V = T S^{-1} T^* = [\ip{f_i}{\tf_j}]_{i,j \in I}$
lies in $\Bc_1(I,a)$.  
Further, by Remark~\ref{selfremark}(b), this occurs if and only if
$(S^{-1/2}(\Fc),a)$ is $\ell^1$-self-localized, where
$S^{-1/2}(\Fc)$ is the canonical Parseval frame.
\end{remark}

\begin{proposition} \label{matrixdecay}
Let $I$ be a countable index set and $a \colon I \to G$ an associated
map such that $D^+(I,a) < \infty$,
and let $K  = \sup_{n \in G} |a^{-1}(n)|$.
Then the following statements hold.

\smallskip
\begin{enumerate}
\item[(a)]
If $V$ has $\ell^1$-decay and $r$ is an associated sequence, then
$V$ maps $\ell^2(I)$ boundedly into itself, with operator norm
$\norm{V} \le K \, \norm{r}_{\ell^1}$.

\medskip
\item[(b)]
The following statements hold:

\smallskip
\begin{enumerate}
\item[i.]
$\Bc_1(I,a)$ is closed under addition and multiplication,

\smallskip
\item[ii.] the following is a norm on $\Bc_1(I,a)$:
$$\norm{V}_{\Bc_1}
\EQ \inf\set{\norm{r}_{\ell^1} : r \text{ is a sequence associated to } V},$$

\smallskip
\item[iii.]
$\Bc_1(I,a)$ is complete with respect to this norm, and

\smallskip
\item[iv.]
we have
\begin{equation} \label{matrixnorm}
\norm{VW}_{\Bc_1} \LE K \, \norm{V}_{\Bc_1} \, \norm{W}_{\Bc_1}.
\end{equation}
\end{enumerate}

\smallskip\noindent
In particular, if $K=1$ then $\Bc_1(I,a)$ is a Banach algebra.

\medskip
\item[(c)]
If $V \in \Bc_1(I,a)$ and $r$ is an associated sequence, then for any
polynomial $p(x) = c_0 + c_1 x + \cdots + c_n x^N$
we have $p(V) \in \Bc_1(I,a)$, and an associated sequence is
$$|c_0| \, \delta + |c_1| \, r + K |c_2| \, (r*r) + \cdots +
K^{n-1} |c_n| \, (r * \cdots * r),$$
where $\delta = (\delta_{0k})_{k \in G}$.

\end{enumerate}
\end{proposition}

\begin{proof}
(a) Given a sequence $c = (c_i)_{i \in I} \in \ell^2(I)$, define
$d \in \ell^2(G)$ by
$$d_n \EQ \sum_{j \in a^{-1}(n)} |c_j|,$$
where we define the sum to be zero if $a^{-1}(n) = \emptyset$.
Note that $\norm{d}_{\ell^2} \le K^{1/2} \, \norm{c}_{\ell^2}$.
Given $i \in I$, we have
\begin{align*}
|(Vc)_i|
\LE \sum_{j \in I} |v_{ij}| \, |c_j|
& \LE \sum_{j \in I} r_{a(i) - a(j)} \, |c_j|
      \\[1 \jot]
& \EQ \sum_{n \in G} \, \sum_{j \in a^{-1}(n)} r_{a(i) - n} \, |c_j|
      \allowdisplaybreaks \\[1 \jot]
& \EQ \sum_{n \in G} r_{a(i) - n} \, d_n
      \\[1 \jot]
& \EQ (r*d)_{a_i}.
\end{align*}
Therefore,
$$\norm{Vc}_{\ell^2}^2
\LE \sum_{i \in I} |(r*d)_{a(i)}|^2
\LE K \, \norm{r*d}_{\ell^2}^2
\LE K \, \norm{r}_{\ell^1}^2 \, \norm{d}_{\ell^2}^2
\LE K^2 \, \norm{r}_{\ell^1}^2 \, \norm{c}_{\ell^2}^2.$$

\medskip
(b) Let $\set{\delta_i}_{i \in I}$ be the standard basis for $\ell^2(I)$.
Suppose $V$, $W \in \Bc_1(I,a)$ with associated sequences $r$, $s$,
and let $c \in \C$.
Then
$$\bigabs{\bigip{(cV+W)\delta_i}{\delta_j}}
\LE |c| \, r_{a(i)-a(j)} + s_{a(i)-a(j)}
\EQ (|c| \, r + s)_{a(i) - a(j)}$$
and
\begin{align*}
|\ip{WV\delta_i}{\delta_j}|
\EQ |\ip{V\delta_i}{W^*\delta_j}|
& \EQ \biggabs{\sum_{k \in I}
      \ip{V\delta_i}{\delta_k} \, \ip{\delta_k}{W^*\delta_j}}
      \\[1 \jot]
& \LE \sum_{k \in I} |\ip{V\delta_i}{\delta_k}| \, |\ip{W\delta_k}{\delta_j}|
      \allowdisplaybreaks \\[1 \jot]
& \LE \sum_{k \in I} r_{a(i) - a(k)} \, s_{a(k) - a(j)}
      \allowdisplaybreaks \\[1 \jot]
& \LE K \, \sum_{n \in G} r_{a(i) - n} \, s_{n - a(j)}
      \\[1 \jot]
& \EQ K \, (r*s)_{a(i) - a(j)}.
\end{align*}
These facts show that $\Bc_1(I,a)$ is an algebra and establish the
norm inequality in \eqref{matrixnorm}.
It is easy to see that $\norm{\cdot}_{\Bc_1}$ is indeed a norm on
$\Bc_1(I,a)$, so it only remains to show that $\Bc_1(I,a)$
is complete with respect to this norm.

Assume that $V_n = [v_{ij}^n]_{i,j \in I}$ for $n \in \N$
forms a Cauchy sequence of matrices in $\Bc_1(I,a)$.
Then, for every $\eps > 0$ there is $N_\eps > 0$ so that for every
$m$, $n \ge N_\eps$ there is a sequence $r^{m,n} \in \ell^1(G)$ such that
$$|v_{ij}^n - v_{ij}^m| \le r^{m,n}_{a(i)-a(j)}
\qquad\text{and}\qquad
\norm{r^{m,n}}_{\ell^1} < \eps.$$
Then for each fixed $i$, $j$, the sequence of entries
$(v_{ij}^n)_{n \in \N}$ is Cauchy, and
hence converges to some finite scalar $v_{ij}$.
Set $V = [v_{ij}]_{i,j \in I}$.

Consider now $\eps_k = \frac{1}{2^k}$ for $k > 0$, and let
$N_k = N_{\eps_k}$ be as above.
Set $N_0 = 0$ and $V^0 = 0$.
Define
$r = \sum_k r^{N_{k+1},N_k}$.
Then $r \in \ell^1(G)$, and
$$|v_{ij}|
\EQ \lim_{k \to \infty} |v^{N_k}_{ij}|
\LE \sum_{k=0}^{\infty} |v^{N_{k+1}}_{ij} - v^{N_k}_{ij}|
\LE r_{a(i)-a(j)}.$$
Hence $V \in \Bc_1(I,a)$,
and it similarly follows that $V^n \to V$ in $\Bc_1(I,a)$.

\medskip
(c) Follows by part~(b) and induction.
\end{proof}

The key to proving Theorem~\ref{selflocthm} is the following fundamental
extension of Wiener's Lemma.
This theorem was proved by Baskakov in \cite{Bas90}
and by Sj\"{o}\-strand in \cite{Sjo95}
(see also \cite{Kur90}, \cite{Bas97}). 

\begin{theorem} \label{sjostrand}
If $V \in \Bc_1(G)$ is an invertible mapping of $\ell^2(G)$ onto itself
then $V^{-1} \in \Bc_1(G)$.
\end{theorem}

\begin{remark}
(a) Sj\"{o}strand proves this result for the case $G = \Z^d$,
but the same technique can be easily applied to the more general groups
we consider in this paper.
Also, Kurbatov proves a more general result for bounded operators
on $\ell^p(\Z^d)$.

\medskip
(b) Theorem~\ref{sjostrand} is similar to Jaffard's Lemma \cite{Jaf90},
which states that if $V$ is invertible on $\ell^2(G)$
and satisfies
$|V_{ij}| \le C \, (1+|i-j|)^{-s}$ for some $C$, $s>0$,
then $V^{-1}$ has the same decay, i.e.,
$|V^{-1}_{ij}| \le C' \, (1+|m-n|)^{-s}$ for some $C'>0$.
Jaffard's Lemma was used by Gr\"{o}chenig
in his development of localized frames in \cite{Gro04}.
~\qed
\end{remark}

Next we define an embedding of the set $\mathbb{F}(I)$ of all frames for $H$
indexed by $I$ into the set $\mathbb{F}(G \times \Z_K)$ of all frames indexed
by $G \times \Z_K$.

\begin{notation} \label{embedding}
Let $I$ be a countable index set and $a \colon I \to G$ and associated
map such that $D^+(I,a) < \infty$.
Let $K  = \sup_{n \in G} |a^{-1}(n)| < \infty$.
For each $n \in G$ let $K_n = |a^{-1}(n)|$, and write
$a^{-1}(n) = \set{i_{nk}}_{k=0}^{K_n-1}$
(it may be the case that $a^{-1}(n)$ is the empty set).
Given a sequence $\Fc = \set{f_i}_{i \in I}$, for each $n \in G$ we set
$$f'_{nk}
\EQ \begin{cases}
    f_{i_{nk}}, & k = 0, \dots, K_n-1, \\
    0,          & k = K_n, \dots, K-1,
    \end{cases}$$
and define $\Fc' = \set{f'_{nk}}_{n \in G, k \in \Z_K}$.
Define $a' \colon G \times \Z_K \to G$ by $a'(i,j) = i$.
\qed
\end{notation}

The following properties are immediate.

\begin{lemma} \label{embedprop}
Let $I$ be a countable index set and $a \colon I \to G$ and associated
map such that $D^+(I,a) < \infty$.
Let $\Fc = \set{f_i}_{i \in I}$ be a frame for $H$.
Then the following statements hold.

\smallskip
\begin{enumerate}
\item[(a)]
$\Fc'$ is a frame for $H$.

\smallskip
\item[(b)]
$(\Fc,a)$ is $\ell^1$-self-localized if and only if
$(\Fc',a')$ is $\ell^1$-self-localized.

\smallskip
\item[(c)]
$(\Fc,a)$ is $\ell^1$-localized with respect to its canonical dual frame
if and only if 
$(\Fc',a')$ is $\ell^1$-localized with respect to its canonical dual frame.

\smallskip
\item[(d)]
If $\tFc$ and $\widetilde{\Fc'}$ denote the canonical duals of
$\Fc$ and $\Fc'$, respectively, then $\widetilde{\Fc'} = (\tFc)'$.
\end{enumerate}
\end{lemma}

Now we can prove Theorem~\ref{selflocthm}.

\begin{proof}[Proof of Theorem \ref{selflocthm}]
By Lemma~\ref{embedprop}, it suffices to consider the case where $I$
is a group of the form given in \eqref{group}, i.e., we can without
loss of generality take $I=G$.
Assume that $\Fc$ is a frame for $H$ such that $(\Fc,a)$
is $\ell^1$-self-localized.
Let $V = [\ip{f_i}{f_j}]_{i,j \in G}$ denote its Gram matrix.
With respect to the algebra $\Bc(\ell^2(G))$ of bounded operators mapping
$\ell^2(G)$ into itself, the spectrum
$\spectrum_{\Bc(\ell^2(G))}(V)$
of $V$ is a closed set contained in
$\set{0} \cup [A,B]$, where $A$, $B$ are the frame bounds of $\Fc$.
On the other hand $V$ belongs to the algebra $\Bc_1(G)$,
and since $\Bc_1(G) \subset \Bc(\ell^2(G))$,
we have the inclusion of spectra
$$\spectrum_{\Bc(\ell^2(G))}(V) \SUBSET \spectrum_{\Bc_1(G)}(V).$$
Theorem~\ref{sjostrand} implies that the converse
inclusion holds true as well, for if
$z \notin \spectrum_{\Bc(\ell^2(G))}(V)$
then $z Id - V$ is an invertible mapping of $\ell^2(G)$ into itself,
and therefore $(z Id - V)^{-1} \in \Bc_1(G)$ by Theorem~\ref{sjostrand}.
Thus
$\spectrum_{\Bc_1(G)}(V)
= \spectrum_{\Bc(\ell^2(G))}(V)
\subset \set{0} \cup [A,B]$.
Let $\Gamma$ denote the circle of radius $B/2$ centered at $(A+B)/2$
in the complex plane.
Then by standard holomorphic calculus \cite{RN90},
the operator
$$V^{\dagger}
\EQ \frac{1}{2\pi i} \int_\Gamma \frac{1}{z} (z Id - V)^{-1} \, dz$$
belongs to $\Bc_1(G)$.
However, the same formula in $\Bc(\ell^2(G))$ defines the pseudoinverse
of~$V$.
Hence $V^{\dagger} \in \Bc_1(G)$, so
$(\tilde{\Fc},a)$ is $\ell^1$-self-localized.
Additionally,
$P_V = V V^\dagger \in \Bc_1(G)$, so
$(\Fc,a)$ is $\ell^1$-localized with respect to its canonical dual
and the associated Parseval frame is $\ell^1$-self-localized.
\end{proof}

\section*{Acknowledgments}
We gratefully acknowledge conversations with Karlheinz Gr\"ochenig and
Massimo Fornasier on localization of frames, and thank them for providing
us with preprints of their articles.
We thank
Hans Feichtinger,
Norbert Kaiblinger,
Gitta Kutyniok,
and
Henry Landau
for conversations regarding the details of our arguments.
We thank Thomas Strohmer and Joachim Toft for bringing the
paper \cite{Sjo95} to our attention,
and Ilya Krishtal for pointing out the paper \cite{Bas97}.

\end{document}